\documentclass{article}
\usepackage{fullpage}
\usepackage{graphics}
 \usepackage{graphicx}	
\usepackage{amssymb}
 \usepackage{amsmath}
 \usepackage{amsthm} 
 \usepackage{lmodern} 
 \usepackage{hyperref}
\usepackage{pdfsync}
\usepackage{yhmath}
\usepackage{amssymb,mathtools}
\usepackage{stmaryrd}
\usepackage{palatino,tipa,graphicx}
\usepackage{lmodern,MnSymbol}

\usepackage{soul}
  \usepackage{xcolor,lmodern}

  \usepackage{calc}
\newsavebox\CBox
\newcommand\hcancel[2][0.5pt]{%
  \ifmmode\sbox\CBox{$#2$}\else\sbox\CBox{#2}\fi%
  \makebox[0pt][l]{\usebox\CBox}%
  \rule[0.5\ht\CBox-#1/2]{\wd\CBox}{#1}}

 \DeclareRobustCommand{\rchi}{{\mathpalette\irchi\relax}}
\newcommand{\irchi}[2]{\raisebox{\depth}{$#1\chi$}} 

 \newcommand{\pd}{\partial} 
 \renewcommand{\H}{\mathcal{H}} 
 
 \newcommand{\Lra}{\Longrightarrow}
 
 \newcommand{\ggm}{\Gamma}
 \newcommand{\R}{{\mathbb{R}}} 
 \newcommand{\8}{\infty}

\renewcommand{\d}{\,{\operatorname{d}}}

\newcommand{\Om}{\Omega}
\newcommand{\TT}{{\mathbb{T}^2}}

\newcommand{\Peri}{\text{Per}}

\renewcommand{\H}{\mathcal{H}} 

\newcommand{\T}{\mathbb{T}}

\newcommand{\vph}{\varphi}
\newcommand{\gr}{\nabla}

 \newcommand{\om}{\Omega}
 
\renewcommand{\vartheta}{\Theta}

 \newtheorem{theorem}{Theorem}[section]
\newtheorem{lemma}[theorem]{Lemma}

\newtheorem{corollary}[theorem]{Corollary}

\DeclareMathOperator{\supp}{supp}

\DeclareMathOperator{\diam}{diam}
	
\newcommand{\vep}{\varepsilon} 

\usepackage{stackengine}
\stackMath
\usepackage{mathtools}
\usepackage{esint}

\title{Decorated phases in triblock copolymers: zeroth- and first-order analysis} 

\author{Stanley Alama
\thanks{Supported by Natural Science and Engineering Research
Council of Canada through the Discovery Grants program. Department of Mathematics and Statistics, McMaster University. E-mail: alama@mcmaster.ca} 
 \qquad Lia Bronsard 
 \thanks{Supported by Natural Science and Engineering Research
Council of Canada through the Discovery Grants program. Department of Mathematics and Statistics, McMaster University. E-mail: bronsard@mcmaster.ca} 
 \qquad Xinyang Lu
 \thanks{Supported by Natural Science and Engineering Research
Council of Canada through the Discovery Grants program. Department of Mathematical Sciences, Lakehead University. Email: xlu8@lakeheadu.ca}
 \qquad Chong Wang
 \thanks{Supported by AMS-Simons Research Enhancement Grant 2024-2027. Department of Mathematics, Washington and Lee University. Email: cwang@wlu.edu}
}

\begin{document}

\date{}
\maketitle

\begin{abstract}

We study a two-dimensional inhibitory ternary system characterized by a free energy functional which combines
an interface short-range interaction energy promoting micro-domain growth with a 
Coulomb-type long-range interaction energy which prevents micro-domains from unlimited spreading.  
Here we consider a scenario in which two species are dominant and one species is
vanishingly small. In this scenario two energy levels are distinguished:  the zeroth-order energy encodes information on
the optimal arrangement of the dominant constituents, while the first-order energy gives the shape of the vanishing constituent.
This first-order energy also shows that,
for any optimal configuration, the vanishing phase must lie on the boundary between the two dominant constituents and form lens clusters also known as vesica piscis.

\end{abstract} 

 \numberwithin{equation}{section}

\section{Introduction}

In this paper, we continue the study started in \cite{ablw, ablw2} of inhibitory ternary systems. Exquisite structures arise in these systems due to competing short-range attractive and long-range repulsive forces. 
An example of such ternary systems are $ABC$ triblock copolymers. An $ABC$ triblock copolymer is a linear-chain molecule consisting of three subchains, joined covalently to each other. A subchain of type $A$ monomer is connected to one of type $B$, which in turn is connected to another subchain of type $C$ monomer. Because of the repulsive forces between different types of monomers, different types of subchain tend to segregate. However, since subchains are chemically bonded in molecules, segregation can lead to a phase separation only at microscopic level, 
where $A, B$ and $C$-rich micro-domains emerge, forming structured patterns \cite{block, zw}.

The energy functional employed in these inhibitory ternary systems is derived from Nakazawa and Ohta's density functional theory for triblock copolymers \cite{microphase, rwtri} in two dimensions. 
We denote by $ \vec{u} = (u_1, u_2, u_0)^{T}$ where each $u_i, i = 0, 1, 2, $ represents the density of type $i$ monomer. In this paper, we pose our problem  
in the flat unit torus $\mathbb{T}^2 =[ - \frac{1}{2}, \frac{1}{2} ]^2$ with periodic boundary conditions. 
The energy functional is defined as
\begin{eqnarray} \label{energyuep}
\mathcal{E}^{\epsilon} (\vec{u})  &:= &   \left [   \sum_{i=0}^2 \frac{ \epsilon }{  M_i }   \int_{\TT}  |\nabla u_i  (\vec{x})|^2  d \vec{x}  + \int_{\TT}  \frac{1}{\epsilon} W(\vec{u})  d \vec{x} \right ]  \nonumber \\
     &&+
  \sum_{i = 0}^2   \sum_{j = 0}^2  \gamma_{ij}  \int_{\TT} \int_{\TT} G_{}( \vec{x}, \vec{y} )\;  (u_i ( \vec{x} ) -  M_i )  \; ( u_j ( \vec{y} ) - M_j )  d \vec{x} d \vec{y},
\end{eqnarray}
which is a diffuse interface model.
Here 
$\epsilon$ is a parameter which is proportional to the thickness of the interfaces, $M_i$ is the area fraction of type $i$ region, 
namely, 
$ \frac{1}{|\TT|} \int_{\TT} u_i ( \vec{x}) d \vec{x} = M_i, i = 0, 1, 2$, 
and
$W$ is a triple-well potential whose three wells represent the three pure phases and are located at $ (1,0,0)^{T},  (0,1,0)^{T}$ and $ (0,0,1)^{T}$.
The long range interaction coefficients $\gamma_{ij}$ form a symmetric matrix $\gamma = [ \gamma_{ij} ] \in \mathbb{R}^{3 \times 3}$
and
$G$ is the Laplace Green's function on $\TT$, of mean zero. In two dimensions, it is known that
\begin{eqnarray*}
G_{} ( \vec{x}, \vec{y} ) = - \frac{1}{2\pi} \log | \vec{x} - \vec{y} | + R ( \vec{x}, \vec{y} ) 
\end{eqnarray*}
for some $R \in C^{\infty} (\mathbb{T}^2)$.
The functional $\mathcal{E}^{\epsilon}$ is defined on the set
 \begin{eqnarray}
\left \{  \vec{u} = (u_1,u_2, u_0)\, : \, u_i \in H^1(\TT; \mathbb{R}) , \ i=0,1,2;  \ \sum_{i=0}^2 u_i = 1; \  \frac{1}{|\TT|} \int_{\TT} u_i ( \vec{x}) d \vec{x} = M_i, i = 1, 2 \right \}, 
\end{eqnarray}
for given $ M_i \in \mathbb{R}, i = 1, 2$.
Using $\Gamma$-convergence \cite{blendCR,miniRW}, the sharp interface limit of $\mathcal{E}^{\epsilon}$ is 
\begin{eqnarray} \label{energyu}
\mathcal{E} ( \Omega_1, \Omega_2, \Omega_0) & :=  & \sum_{0 \leq i < j \leq 2} \sigma_{ij}   \mathcal{H}^1 ( \partial \Omega_i \cap  \partial \Omega_j) \nonumber \\
 && + \sum_{i = 0}^2   \sum_{j = 0}^2   \gamma_{ij}  \int_{\TT} \int_{\TT} G_{}(\vec{x}, \vec{y})\;  (\rchi_{\Om_i} (\vec{x}) - M_i)  (\rchi_{\Om_j} (\vec{y}) - M_j ) \;   d\vec{x} d\vec{y},
\end{eqnarray}
where 
the densities $u_i$ are replaced by phase domains described by characteristic functions $\rchi_{\Om_i}$ of sets $\Om_i, \ i=0,1,2,$ which partition $\TT$ 
into disjoint sets.

The constants $\sigma_{ij}$, $i\neq j$, represent {\it surface tension} along the interfaces separating the phase domains, and they are positive, material-dependent constants. 
In this paper we discuss a special case $\sigma_{01} = \sigma_{02} = \sigma_{12}= 1$, which arises (for instance) when $W$ is symmetric with respect to permutations of the $e_i$-axes.

 \begin{figure}[!htb]
\centering
  \includegraphics[width=4cm]{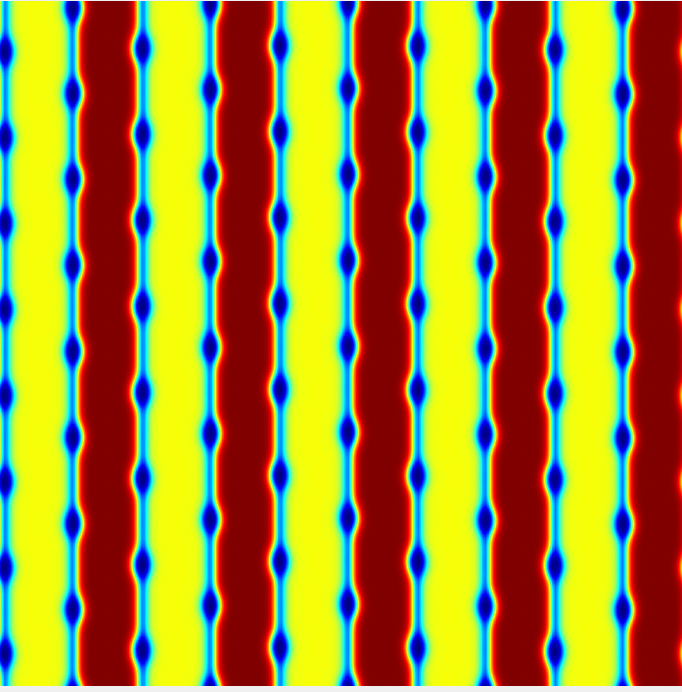} 
  \caption{The Numerical Simulation: decorated phases of ABC triblock copolymers. Type A micro-domains are in red, type B are in yellow, and type C are in blue.}
   \label{numeric1}
\end{figure}

As in \cite{ablw, ablw2}, we explore an asymptotic regime of ternary systems. However, unlike the configuration with one dominant constituent and two vanishing ones in \cite{ablw, ablw2}, we focus on a scenario where two constituents are dominant while one diminishes; see Figure \ref{numeric1}.
 We study the ``droplet'' scaling and introduce a new parameter $\eta$ which is to represent the characteristic length scale of the droplet components as in \cite{bi1, ablw, ablw2, ablw3}.
Thus, areas scale as $\eta^2$, and so we choose area constraints on $(\Omega_1, \Omega_2,\Omega_0)$ as follows
\begin{eqnarray}
  |\Omega_1| =  |\Omega_2| = \frac{1-\eta^2M}{2}, \qquad |\Omega_0|= \eta^2M,
  \label{choice of masses of backgrounds}
\end{eqnarray}
for some fixed $M$. This choice can be replaced by others. For example,
\begin{eqnarray}
  |\Omega_i| =  a_i  (1 - \eta^2 M ) , \;  i = 1, 2, \quad  a_1 + a_2 = 1 , \qquad |\Omega_0|= \eta^2M .
\end{eqnarray}
For the choice of masses in \eqref{choice of masses of backgrounds}, let
\begin{eqnarray}
E_\eta(\om_1,\om_2,\om_0) := \mathcal{E} (\om_1,\om_2,\om_0) .
\end{eqnarray}

We remark that the main differences between this paper and \cite{bi1} are that when a minority phase cluster lies on the boundary between the majority backgrounds, it removes a piece of perimeter from the latter, therefore, as proven in  \cite{alama2023lens}, the optimal shape is a lens instead of the sphere.


\subsection*{Heuristics of scaling}

	We need first to get the orders of 
	the interaction coefficients, $\gamma_{ij}$, so that the local and nonlocal term interact at the same order. 
	To this aim, it is helpful to think of $\om_1$, $\om_2$
	as ``rectangles'' with side lengths of size $O(1)$, and $\om_0$ as small droplets
	(or balls) of mass $O(\eta^2)$ (and hence radius $O({\eta})$). 
	Note all the terms in the sum
	\begin{align*}
	E_L(\om_1,\om_2) &:=  \H^1( \pd \om_1 \cap \pd \om_2 )
	+ \sum_{i,j =1}^2  \gamma_{ij} 
	\int_{\TT} \int_{\TT} G_{}(\vec{x}, \vec{y})\;  (\rchi_{\Om_i} (\vec{x}) - M_i )  (\rchi_{\Om_j} (\vec{y}) - M_j ) \;   d\vec{x} d\vec{y},
	\end{align*}
	are of order $O(1)$. Therefore, we impose
	\begin{align}
	\label{zeroeth order parameters}
 \gamma_{11}, \gamma_{12}, \gamma_{22} =O(1). 
	\end{align}
	Next, similar to the calculations in [\cite{bi1}, (3.6)] we note that, it is possible to choose $\Omega_0$ such that 
	\begin{align*}
	\H^1( \pd \om_1 \cap \pd \om_0 ),\  \H^1( \pd \om_2 \cap \pd \om_0 ) &=O({\eta}),\\
	\int_{\TT} \int_{\TT} G_{}(\vec{x}, \vec{y})\;  (\rchi_{\Om_i} (\vec{x}) - M_i )  (\rchi_{\Om_0} (\vec{y}) - M_0 ) \;   d\vec{x} d\vec{y} &=O(\eta^2),\qquad i=1,2, \\
		\int_{\TT} \int_{\TT} G_{}(\vec{x}, \vec{y})\;  (\rchi_{\Om_0} (\vec{x}) - M_0 )  (\rchi_{\Om_0} (\vec{y}) - M_0 ) \;   d\vec{x} d\vec{y}  &=O(\eta^4|\log \eta|). 
	\end{align*}
Intuitively, in the order $O(\eta^4|\log \eta|)$, $O(\eta^4)$ comes from the fact that $\Omega_0$ has mass $O( \eta^2) $, while the $ O ( |\log \eta | )$ comes from the logarithmic interaction, combined with
the fact that diameter of $\Om_0$ is the order of $O (\eta)$.

Therefore, in order for the interaction terms involving $\Omega_0$
 to have the same order as its perimeter, we impose
	\begin{align}
	\label{first order term}
\gamma_{i0} = \frac{\Gamma_{i0}}{\eta}, \quad i = 1,2,
	\end{align}
and
	\begin{align}
	\label{self interaction type two}
	\gamma_{00} = \frac{\Gamma_{00}}{\eta^3|\log\eta|}.
	\end{align}
Since $G$ is of mean zero, we have
\begin{eqnarray*}
\int_{\TT} \int_{\TT} G_{}(\vec{x}, \vec{y})\;  (\rchi_{\Om_i} (\vec{x}) -  M_i )  (\rchi_{\Om_j} (\vec{y}) - M_j ) \;   d\vec{x} d\vec{y} = 
 \int_{\om_i} \int_{\om_j} G_{}(\vec{x}, \vec{y}) d\vec{x} d\vec{y}.
\end{eqnarray*}
Thus we can write the energy $E_\eta$ as 
	\begin{align}
	E_\eta & (\om_1,\om_2,\om_0)  =  \underbrace{ E_L(\om_1,\om_2)}_{=O(1)} \notag\\
	&+ \eta\bigg[ \underbrace{
	\frac{1}{\eta}\sum_{i=1}^2 \H^1( \pd \om_i \cap \pd \om_0 )	
		+ \sum_{i=1}^2  \frac{ \Gamma_{i0}}{\eta^2} \int_{\om_i} \int_{\om_0} G_{}(\vec{x}, \vec{y}) d\vec{x} d\vec{y}
	+\frac{\Gamma_{00}}{\eta^4 |\log \eta|} \int_{\om_0} \int_{\om_0} G_{}(\vec{x}, \vec{y}) d\vec{x} d\vec{y} }_{=O(1)} \bigg].
	\label{energy with correct scaling}
	\end{align}
       We remark that,
	when the parameters in \eqref{zeroeth order parameters} are sufficiently small,
	the global minima of $E_L$ are lamellars (see \cite{st}). Here,
	 the choices
	on the masses of the backgrounds
	in 	
	\eqref{choice of masses of backgrounds},
	 and on the interaction coefficients 
in \eqref{zeroeth order parameters}, 
can actually be relaxed: the main goal is 
to satisfy the hypotheses 
of \cite[Proposition~2.1]{st}, that is, the two majority phases, i.e. $\Omega_1$, $\Omega_2$, are a local minima for $E_L$,
which provides $C^{3,\alpha}$ (for some $\alpha\in (0,1)$)
regularity for the (reduced) boundary of $\Omega_1$ and $ \Omega_2$
 for minima. This
ensures that, under blow-up, these boundaries become straight lines,
which
 will be crucial in the optimality of lenses (see Subsection \ref{optimal lenses}
 below).

In previous ternary studies \cite{double, doubleAs, stationary, disc, wrz, ablw, ablw2}, one species was the background, and the other two interacted at the same order, making it possible to set the former to be the ``ground level", thus essentially removing it from the interaction altogether. In the current paper, however, the presence of two majority phases, plus a minority one, introduces interactions at different levels. Thus the reduction of the interaction coefficient matrix to a two by two matrix is not possible, and we will keep the original three by three interaction coefficient matrix $\gamma$.
This matrix was first derived by Nakazawa and Ohta using mean field theory [\cite{microphase}, Equations (2.23) and (A.7)] and later by Ren and Wei [\cite{rwtri}, Equation (4.21)]. Unlike the reduced matrix, the 
three by three matrix includes negative $\gamma_{ij}$
values. Therefore, in this paper, we extend the analysis to include cases where the interaction coefficients are negative, provided they satisfy a smallness condition.


The paper is organized as follows: In Section~\ref{zeroth}, we prove a zeroth-order $\Gamma$-convergence result (Lemma \ref{zeroth order gamma limit}). In Section ~\ref{first}, we discuss the limit function of the first-order
$\Gamma$-convergence
(see \eqref{1st order Gamma limit energy} below), the optimal shape of the minority phase
(the so-called ``optimal lenses'', see
	 	\cite{alama2023lens}), and then we provide the proof of the first-order $\Gamma$-convergence
(Theorem \ref{1st order Gamma convergence}).
In Section~\ref{uniform}, we show that in optimal configurations, the minority phase is spread over a finite number of equally sized masses
(Theorem \ref{uniform lenses theorem}).
In Section 5, we discuss the case of negative interaction coefficients.

 Due to its complexity, the
mathematical study of \eqref{energyu} is still in its early stages. One-dimensional stationary points to the Euler-Lagrange equations of \eqref{energyu} were found in \cite{lameRW, blendCR}.
Two and three dimensional stationary configurations were recently studied in \cite{double, doubleAs, stationary, disc, evolutionTer}, and global minimizers in \cite{ablw, ablw2}.

While the mathematical interest in triblock copolymers via the energy functional \eqref{energyu} is relatively recent, there has been much progress in the mathematical analysis of nonlocal binary systems.  Much early work concentrated on the diffuse interface Ohta-Kawasaki density functional theory for diblock copolymers \cite{equilibrium, nishiura, onDerivation},
\begin{eqnarray} \label{energyB}
\mathcal{E}^{} (\Omega) :=  
 \H^1( \pd \om) + \gamma \int_{\T^n} \int_{\T^n} G_{}(\vec{x}, \vec{y})\;  (\rchi_{\Om} (\vec{x}) - M)  (\rchi_{\Om} (\vec{y}) - M) \;   d\vec{x} d\vec{y}
\end{eqnarray}
with a single mass or volume constraint. The dynamics for a gradient flow for \eqref{energyB} with small volume fraction were developed in \cite{hnr, gc}. 
All stationary solutions to the Euler-Lagrange equation of \eqref{energyB} in one dimension were known to be local minimizers \cite{miniRW}, and
many stationary points in two and three dimensions have been found that match the morphological phases in diblock copolymers \cite{oshita, many, spherical, oval, ihsan, Julin3, cristoferi, afjm}.
The sharp interface nonlocal isoperimetric problems have been the object of great interest, both for applications and for their connection to problems of minimal or constant curvature surfaces \cite{otto, bi1, st,  knupfer1, knupfer2, Julin, ms}.
Global minimizers of \eqref{energyB}, and the related Gamow's Liquid Drop model describing atomic nuclei, were studied in \cite{ muratov, bi1, st, GMSdensity, fl} for various parameter ranges.
Variants of the Gamow's liquid drop model with background potential or with an anisotropic surface energy replacing the perimeter, are studied in \cite{ABCT1,luotto, cnt}.
Higher dimensions are considered in \cite{BC, cisp}.
Applications of the second variation of \eqref{energyB} and its connections to minimality and $\Gamma$-convergence are to be found in \cite{cs,afm,Julin2}.
 The bifurcation from spherical, cylindrical and lamellar shapes with Yukawa instead of Coulomb interaction has been done in \cite{fall}.
Blends of diblock copolymers and nanoparticles \cite{nano, ABCT2} and blends of diblock copolymers and homopolymers are also studied by \cite{BK,blendCR}.
Extension of the local perimeter term to nonlocal $s$-perimeters is studied in \cite{figalli}.


\vspace{0.2cm}

\section{Zeroth-order Gamma convergence} \label{zeroth}

	In this section, we show the Gamma-convergence of $E_\eta$ to $E_L$
	as $\eta\to 0$.
	
	
	\begin{lemma}
		\label{zeroth order gamma limit}
		It holds $E_\eta \overset{\ggm}{\to}E_L$. 
	\end{lemma}

\begin{proof}
		We need to prove compactness, $\ggm-\liminf$, and $\ggm-\limsup$ inequalities.
		
		\medskip
		
		\textbf{Step 1. Compactness.}
		 Consider a sequence $(\om_{1,\eta},\om_{2,\eta},\om_{0,\eta})$,
		such that 
		$$ |\om_{1,\eta}|=|\om_{2,\eta}|=\frac{1-\eta^2 M}{2},\quad
		|\om_{0,\eta}|=\eta^2 M,
		\qquad \sup_\eta E_\eta(\om_{1,\eta},\om_{2,\eta},\om_{0,\eta})=:C <+\8.$$ 
		It's straightforward to check 
		that
		\begin{align}
				&\sum_{i=1}^2 \frac{\ggm_{i0}}{\eta } \int_{\om_{i,\eta}} \int_{\om_{0,\eta}} G_{}(\vec{x}, \vec{y}) d\vec{x} d\vec{y}
		\notag\\
		&
		\ge \sum_{i=1}^2 \frac{\ggm_{i0}}{\eta } |\om_{i,\eta}||\om_{0,\eta}|
		\left[-\frac{1}{2\pi} \log \diam \T^2 +  \underset{\T^2}{\inf}   R\right] 
		= (\ggm_{10} + \ggm_{20}) \eta M \frac{1-\eta^2 M }{2}
		\left[-\frac{1}{2\pi} \log \diam \T^2 +  \underset{\T^2}{\inf} R  \right] ,
		\label{01-2 interaction is negligible}\\
		\text{ and }   \nonumber \\
		&\frac{\ggm_{00}}{\eta^3 |\log \eta|}\int_{\om_{0,\eta}} \int_{ \om_{0,\eta}} G_{}(\vec{x}, \vec{y}) d\vec{x} d\vec{y}
		\notag\\
		&
		\ge \frac{\ggm_{00}}{\eta^3 |\log \eta|}|\om_{0,\eta}|^2
		\left[-\frac{1}{2\pi} \log \diam \T^2 + \underset{\T^2}{\inf}  R\right] 
		=\frac{\ggm_{00} \eta M^2 }{|\log \eta|}
		\left[-\frac{1}{2\pi} \log \diam \T^2 +  \underset{\T^2}{\inf}  R\right] ,
		\label{22 interaction is negligible}
		\end{align}
		and hence 
		\begin{align*}
				 \lim_{\eta\to 0} 
				 \sum_{i=1}^2 \frac{\ggm_{i0}}{\eta } \int_{\om_{i,\eta}} \int_{\om_{0,\eta}} G_{}(\vec{x}, \vec{y}) d\vec{x} d\vec{y}
		  &\ge 0,\\
				 \lim_{\eta\to 0} 
	\frac{\ggm_{00}}{\eta^3 |\log \eta|} \int_{\om_{0,\eta}} \int_{ \om_{0,\eta}} G_{}(\vec{x}, \vec{y}) d\vec{x} d\vec{y}		
			&\ge 0. 
		\end{align*}
		Since $\sum_{i=1}^2 \H^1( \pd \om_i \cap \pd \om_0 )	\ge 0$,
		condition  $\sup_\eta E_\eta(\om_{1,\eta},\om_{2,\eta},\om_{0,\eta})<+\8$ now implies 
		$E_L(\om_{1,\eta}, \om_{2,\eta})<+\8$.
		Again, the interaction term satisfies
		\begin{align}
		\sum_{i,j =1}^2  \gamma_{ij} \int_{\om_{i,\eta}} \int_{ \om_{j,\eta}} G_{}(\vec{x}, \vec{y}) d\vec{x} d\vec{y}	
		\ge \sum_{i,j =1}^2 \gamma_{ij} |\om_{i,\eta}|| \om_{j,\eta}| \left[-\frac{1}{2\pi} \log \diam \T^2 +  \underset{\T^2}{\inf}  R\right]
		>-\8. \label{interaction term bounded from below}
		\end{align}
		Hence letting $ u_{i,\eta} = \chi_{\Omega_{i, \eta}}$, 
		this implies
		\[\sum_{i,j=1 }^2 \int_{\T^2} |\gr u_{i,\eta}| \d \vec{x}<+\8.\]	
		As the total masses $| \om_{i,\eta}|= || u_{i,\eta} ||_{L^1(\T^2)}$, $i=1,2 $,
		are bounded, we get that $u_{i,\eta}$ are uniformly bounded in $BV(\T^2)$, so they
		converge to limit functions $u_{i} = \chi_{\Omega_i} $ in the weak topology of $BV(\T^2)$. It also follows
		that $\sum_{i=0}^2 u_{i}=\sum_{i=0}^2 u_{i,\eta}=1$ and $|\om_i\cap \om_j|=0$, $i\neq j$.

		\medskip
		
		\textbf{Step 2. $\ggm-\liminf$ inequality.} We need to show that given $u_{i,\eta} = \chi_{\Omega_{i, \eta}}  $
		converging to $u_{i} =\chi_{\om_i} $ in the weak topology of $BV(\T^2)$, it holds
		\begin{align}
		\label{gamma liminf inequality}
		\liminf_{\eta\to 0} E_\eta(\om_{1,\eta},\om_{2,\eta},\om_{0,\eta}) \ge 
		E_L(\om_{1},\om_{2}) . \qquad 
		\end{align}
		By \eqref{01-2 interaction is negligible}, \eqref{22 interaction is negligible}, and the fact that $\sum_{i=1}^2 \H^1( \pd \om_i \cap \pd \om_0 )\ge 0$, we get
		\[	\liminf_{\eta\to 0} E_\eta(\om_{1,\eta},\om_{2,\eta},\om_{0,\eta}) \ge 
		\liminf_{\eta\to 0} E_L(\om_{1,\eta},\om_{2,\eta}),\]
		so it suffices to show
		\[\liminf_{\eta\to 0} E_L(\om_{1,\eta},\om_{2,\eta}) \ge E_L(\om_{1},\om_{2}).\]
		This is true since
		\[ E_L(\om_{1,\eta},\om_{2,\eta})=
		\sum_{i=1 }^2\int_{\T^2} |\gr u_{i,\eta}| \d \vec{x}+
		\sum_{i,j=1 }^2  \gamma_{ij} \int_{\om_{i,\eta}} \int_{ \om_{j,\eta}} G_{}(\vec{x}, \vec{y}) d\vec{x} d\vec{y} , \]
		where the first term is clearly lower semicontinuous with respect to the weak convergence
		in $BV(\T^2)$, and the second term
		\begin{align*}
		\lim_{\eta\to 0}  \int_{\om_{i,\eta}} \int_{ \om_{j,\eta}} G_{}(\vec{x}, \vec{y}) d\vec{x} d\vec{y}  &=
		\lim_{\eta\to 0}
		\int_{\T^2} \int_{ \T^2} G_{}(\vec{x}, \vec{y})  u_{i,\eta}(\vec{x})u_{j,\eta}(\vec{y}) \d \vec{x} \d \vec{y} \\
		&=\lim_{\eta\to 0}\int_{\T^2} u_{i,\eta}( \vec{x})  \bigg[\int_{\T^2}G(\vec{x}, \vec{y} ) u_{j,\eta}(\vec{y})  \d \vec{y} \bigg] \d \vec{x}\\
		&
			=\lim_{\eta\to 0}\int_{\T^2} u_{i,\eta}(\vec{x})  [G(\vec{x}, \cdot)* u_{j,\eta}]  \d \vec{x}\\
		&= \int_{\T^2} u_{i}(\vec{x})  [G(\vec{x}, \cdot)* u_{j}] \d \vec{x}
		=\int_{\om_{i}} \int_{ \om_{j}} G(\vec{x}, \vec{y})\d \vec{x} \d \vec{y} 
		\end{align*}
		since $u_{i,\eta} \to u_{i}$ and $G(\vec{x}, \cdot)* u_{j,\eta}\to G(\vec{x}, \cdot)* u_{j}$
		strongly in $L^2(\T^2)$.
		
		\medskip
		
		
		\textbf{Step 3. $\ggm-\limsup$ inequality.}
		We need to show that given $\om_{1},\om_{2}$, there exists a sequence $u_{i,\eta}$
		converging to $u_{i}=\chi_{\om_i}$ in the weak topology of $BV(\T^2)$, such that
		\begin{align}
                  ||u_{1,\eta}||_{L^1(\T^2)}= ||u_{2,\eta} ||_{L^1(\T^2)}&=\frac{1-\eta^2 M}{2},\quad
		||u_{0,\eta}||_{L^1(\T^2)}=\eta^2 M,\notag\\
		\limsup_{\eta\to 0} E_\eta(\om_{1,\eta},\om_{2,\eta},\om_{0,\eta}) &\le 
		E_L(\om_{1},\om_{2}).		\label{gamma limsup inequality}
		\end{align}
		Assume $E_L(\om_{1},\om_{2})<+\8$, as otherwise the inequality is trivial.
		We construct $u_{i,\eta}$ by perturbing $u_{i}$: 
		\begin{enumerate}
			\item choose $x_i$ to be a density point for
			$\om_i$, that is
			\[ \lim_{\vep\to 0}\frac{|B(x_i,\vep) \cap \om_i|}{|B(x_i,\vep)|} =1 ,\qquad i=1,2. \]
			Here $\varepsilon$ is a parameter that, as we will show later, can be chosen to
			vanish as $\eta\to 0$.
			
			\item Replace $B(x_i,\vep) $ with a core-shell (i.e. annular region) where the inner circle
			$B(x_i,\vep_i)$ is filled with type $i$ constituent, with $\vep_i$ chosen such that 
			$|B(x_i,\vep_i)|=|B(x_i,\vep) \cap \om_i|$. Then fill the outer shell $B(x_i,\vep)\setminus
			B(x_i,\vep_i)$ with the other type constituent. Note this construction 
			requires only $\vep< |x_1-x_2|/2 $, and does not 
			alter
			the total masses of each constituent type in each $B(x_i,\vep)$.
			
			\item Fill the balls $B(x_i, \eta\sqrt{M/2\pi} ) $ with type 0 constituent. 
			Without loss of generality, we can choose 
			$\eta$ sufficiently small 
			$\eta\sqrt{M/2\pi}<\min \{\vep_1,\vep_2\}$.
			This construction
			adds mass $\eta^2M$ of type 0 constituent, and removes mass $\eta^2M/2$ of both type
			I and II constituent. The resulting configuration then satisfies the mass constraints
			\[ |\om_{1,\eta}|=|\om_{2,\eta}|=\frac{1-\eta^2M}{2},\qquad |\om_{0,\eta}|=\eta^2M. \]
		\end{enumerate}
		This construction can create, at most the following six perimeters: 
		\[\pd B(x_i,\vep_i),\qquad
		\pd B(x_i,\vep),\qquad \pd B(x_i,\eta\sqrt{M/2\pi}),\qquad i=1,2.\]
		Thus the total perimeter is increased by at most $+ 8\pi \vep + 2\eta\sqrt{2M\pi}$.
		Note the only requirement on $\vep $ is $\eta\sqrt{M/2\pi}<\min \{\vep_1, \vep_2 \}
		\le \max \{\vep_1, \vep_2 \}\le \varepsilon$,
		and since $x_i$ are density points for $\om_i$, for all sufficiently small $\vep$
		we have $\vep_i\ge \vep/2$. Thus we only need to choose $\vep\ge \eta\sqrt{2M/\pi}$, hence
		$\vep$ can be chosen arbitrarily small as $\eta\to 0$. Therefore, the
		maximum increase in perimeter is vanishing
		as  $ \eta\to 0$.
		Finally, note that by construction
		$u_{i,\eta} = \chi_{\om_{i,\eta} } \to u_i = \chi_{\om_i}$ strongly in $L^p(\T^2)$
		for any $p<+\8$, and
		it follows that the interaction parts converge too.
		
	\end{proof}

\vspace{0.2cm}


\section{First-order Gamma convergence} \label{first}
		
		We now analyze the $\Gamma$-convergence of $\frac{E_\eta-E_L}{\eta}$.
		Qualitatively, this includes the terms accounting for the perimeter of type 0 constituent regions, their interaction with the two backgrounds $\Omega_1$ and $\Omega_2$, and their self interactions.

\subsection{Qualitatively observations}
		
		In this section we assume that
			 $\gamma_{10}=\gamma_{20}$, namely that the interaction between type I and 0 constituents is the the same as that between type II and 0 constituents. 
			 This means that, for the first order $\Gamma$-convergence,
			the interactions between the minority phase and the two majority constituents are indistinguishable.

			\medskip
			
			To simplify the presentation,
			 we consider the case of a {\em finite} partition since such sets can approximate the general case of countable partition arbitrarily well: for any configuration $\Omega $ with countably many components $\Omega_{i,k}$, $i = 0, 1, 2$,
			$k \ge 1$, the energy of the sequence of (truncated) sets 
			\[ \frac{ |\Omega |}{ |  \bigcup_{k=1}^N \Omega_{i, k } | }  \bigcup_{k=1}^N \Omega_{i, k},  \]
			rescaled to have the same mass of $\Omega$, converges to the energy of $\Omega$. So let
			\[ \Omega_{0,\eta} = \bigsqcupdot_{k=1}^n\Omega_{0,\eta,k} ,\qquad
			|\Omega_{0,\eta,k} |=\eta^2 m_k,\qquad \sum_{k=1}^n m_k=M,\]
			where $\bigsqcupdot$ represents the disjoint union and each $\Omega_{0,\eta,k}$ is a connected component,
			satisfying
			\begin{eqnarray*}
			\label{positive minimum distance between different clusters}
						\min_{k\neq l }\{ |\vec{x}- \vec{y}|:  \vec{x} \in \Omega_{0,\eta,k} , \, \vec{y} \in \Omega_{0,\eta,l}  \}\ge d_0\gg \eta.
			\end{eqnarray*}
			%
			A priori, we do not know the shape of each $\Omega_{0,\eta,k}$. However,
			we can always replace $\Omega_{0,\eta,k}$ with a ball $B$ of identical mass,
			which would yield an upper bound on its energy contribution to the overall energy as follows: 	 
			let $\Gamma := \Gamma_{10} = \Gamma_{20}$ (defined in \eqref{first order term}), 
			\begin{align}
			\H^1(\partial B)&
			+\Gamma\bigg[  \frac{1}{\eta}      
			\int_{B}  \int_{\mathbb{T}^2 \setminus 
			(\sqcupdot_{l=1, l \neq k}^n\Omega_{0,\eta,l} \cup B )}
			 G(\vec{x}, \vec{y} )\d  \vec{x} \d \vec{y}  \notag \\
		&\qquad+
		\frac{1}{\eta^3 |\log \eta|}\int_{B} \int_{ \sqcupdot_{l=1, l\neq k}^n\Omega_{0,\eta,l}} G(\vec{x}, \vec{y})\d \vec{x} \d \vec{y} 
		+
		\frac{1}{\eta^3 |\log \eta|}\int_{B} \int_{B} G(\vec{x}, \vec{y})\d \vec{x} \d \vec{y}  \bigg] \notag\\
		&=
		2 \eta\sqrt{\pi m_k} + \Gamma\bigg[\frac{1}{\eta}\int_{B} \int_{\mathbb{T}^2 \setminus
		B } 
		G(\vec{x}, \vec{y})\d \vec{x} \d \vec{y}   \notag\\
		&\qquad  +\Big(\frac{1}{\eta^3 |\log \eta|}-
		\frac{1}{\eta}\Big)\sum_{l=1, l\neq k}^n\int_{B} \int_{ 
			\Omega_{0,\eta, l}} 
			G(\vec{x}, \vec{y})\d \vec{x} \d \vec{y} \notag\\
		&\qquad
		+
		\frac{1}{\eta^3 |\log \eta|}\int_{B} \int_{ B
		} G(\vec{x}, \vec{y})\d \vec{x} \d \vec{y}  \bigg] 
	\notag\\
	\end{align}
	\begin{align}
	&=
	2 \eta\sqrt{\pi m_k} + \Gamma\bigg[\frac{1}{\eta}\int_{B} \int_{ \mathbb{T}^2} G(\vec{x}, \vec{y})\d \vec{x} \d \vec{y} -
	\frac{1}{\eta}\int_{B}  
		\int_{ B } G(\vec{x}, \vec{y})\d \vec{x} \d \vec{y} \notag\\
	&\qquad  +\Big(\frac{1}{\eta^3 |\log \eta|}-
	\frac{1}{\eta}\Big)\sum_{l =1, l\neq k}^n\int_{B}  
		\int_{ \Omega_{0,\eta, l}} G(\vec{x}, \vec{y})\d \vec{x} \d \vec{y}  \notag\\
	&\qquad
	+
	\frac{1}{\eta^3 |\log \eta|}\int_{B} \int_{ B
	} G(\vec{x}, \vec{y})\d \vec{x} \d \vec{y}  \bigg].
	\label{energy contribution of a ball - decomposition}
			\end{align} 
			 %
			 In view of \eqref{positive minimum distance between different clusters}, similar to the calculations in [\cite{bi1}, (3.6)], we have
			 \begin{align}
			 \frac{1}{\eta}\int_{B} \int_{
			 	B } G(\vec{x}, \vec{y})\d \vec{x} \d \vec{y} 
		 	&=O( \eta^3|\log  \eta| ),
		 	\label{negligible term 1}\\
			 			 \Big(\frac{1}{\eta^3 |\log \eta|}-\frac{1}{\eta}\Big)\sum_{l =1, l \neq k}^n\int_{B} \int_{ 
			 	\Omega_{0,\eta,l}} G(\vec{x}, \vec{y})\d \vec{x} \d \vec{y} &= O\Big(\frac {\eta}{|\log \eta|}\Big) .
			 \label{negligible term 2}
			 \end{align}
			 Moreover, since $m_k\le M$, hence
			 $m_k\le \sqrt{M m_k}$, the other three terms in \eqref{energy contribution of a ball - decomposition}, namely
			 \[ 2 \eta\sqrt{\pi m_k}, \quad \frac{1}{\eta}\int_{B} \int_{ \mathbb{T}^2} G(\vec{x}, \vec{y})\d \vec{x} \d \vec{y}, \quad \frac{1}{\eta^3 |\log \eta|}\int_{B} \int_{ B
	} G(\vec{x}, \vec{y})\d \vec{x} \d \vec{y} \]
			  are all of order $O(\eta)$ as $B$ has mass $O(\eta^2)$ and the $ O ( |\log \eta | )$ is cancelled due to the logarithmic interaction, combined with
the fact that the diameter of $B$ is the order of $O (\eta)$, thus \eqref{negligible term 1}
			 and \eqref{negligible term 2}
			 	are negligible in the sum \eqref{energy contribution of a ball - decomposition}.
			 Consequently, there exists a computable constant
			 $C_{E}(\Gamma,M)$ such that the energy contribution of $B$ \eqref{energy contribution of a ball - decomposition}
			 is bounded from above by
			 \begin{align}
			 \label{energy contribution of the ball}
			  C_E(\Gamma,M)\eta \sqrt{m_k}.
			 \end{align}
			 As a consequence, if the original configuration were optimal, then the energy
			 contribution of $\Omega_{0,\eta,k}$ cannot exceed \eqref{energy contribution of the ball}.
			In conclusion, such point bound on the energy contribution of each $\Omega_{0,\eta,k}$ yields
			a bound on
			their perimeter too, i.e. there exists another constant
			$C_{P}(\Gamma,M)$ such that
			\begin{align}
			\label{bounded perimeter}
						\H^1(\partial \Omega_{0,\eta,k}) \le C_{P}(\Gamma,M)\eta \sqrt{m_k}, \qquad \text{for all }k.
			\end{align}

\subsection{Expected form of the 1st order $\Gamma$-limit energy}

		We now determine the form of the limit function. Similarly to what was done in \cite{bi1,ablw}, we do a blow up of order $1/\eta$:
		 set
		\[v_{0,\eta}:\mathbb{T}^2 \longrightarrow \mathbb{R},\qquad   v_{0,\eta}:=\frac{1}{\eta^2} \sum_{k} \chi_{\Omega_{0,\eta,k}} ,\]
		which is the analogue of \cite[(1.6)]{bi1},
		hence 
		\[ \|v_{0,\eta}\|_{L^1(\mathbb{T}^2)} =M,\qquad \H^1(\partial \Omega_0) = \int_{\mathbb{T}^2} |\nabla v_{0,\eta}| \text{d} \vec{x} = \sum_k \int_{\mathbb{T}^2} |\nabla v_{0,\eta,k}| \text{d} \vec{x}. \]
		As $\eta\to 0$, each $\frac{1}{\eta^2} \chi_{\Omega_{0,\eta,k}}$ should concentrate
		to a Dirac mass, hence the $v_{0,\eta}$ should converge, in the weak topology of the space of Radon measures on $\mathbb{T}^2$, to a sum of Dirac masses.
		Define also 
		\[z_{0,\eta}: 
		\Big[ -\frac{1}{2\eta}, \frac{1}{2\eta}\Big]^2\subseteq \mathbb{R}^2\longrightarrow \mathbb{R},\quad \text{via}  \quad z_{0,\eta}(x):= \sum_{k} \rchi_{\Omega_{0,\eta,k}}(\eta x
		),
		\]
		which is the analogue of \cite[(3.2)]{bi1}.
		Next observe that there is a close similarity between 
		\[\sum_{i=1}^2 \H^1( \pd \om_i \cap \pd \om_0 )	
		+\frac{\Gamma_{00}}{\eta^4 |\log \eta|}\int_{\om_0} \int_{ \om_0} G(\vec{x}, \vec{y} )\d  \vec{x} \d \vec{y} \]
		and $E^{2d}_\eta$ defined in \cite[(1.6)]{bi1}.
		Thus following their approach we define
		\begin{align}
		\bar{e}(m):=\inf \Big\{ \sum_k e(m_k):m_k\ge 0 ,\sum_k m_k =M \Big\},\qquad
		e(m):=c_1 \sqrt{m} + c_2m^2,\label{limit functions e}
		\end{align}
		for some constants $c_1,c_2>0$,
		which is the analogue of \cite[(6.1)]{bi1}.
		On the other hand, the first order energy $\frac{E_\eta-E_L}{\eta}$
		has also the term
		\begin{align}
		\label{energy - 1st order}
				\sum_{i=1}^2 \frac{\Gamma_{i0}}{\eta^2} \int_{\om_i} \int_{ \om_0} G(\vec{x}, \vec{y} )\d  \vec{x} \d \vec{y},
		\end{align}
		which can be rewritten as
		\begin{align*}
		\frac{\Gamma}{\eta^2} & \sum_{i=1}^2 \int_{\om_i} \int_{ \om_0} G(\vec{x}, \vec{y} )\d  \vec{x} \d \vec{y}
		=
		\frac{\Gamma}{\eta^2}\int_{\om_0} \int_{\mathbb{T}^2 \setminus
		\Omega_0}  G(\vec{x}, \vec{y} )\d  \vec{x} \d \vec{y} 
	=
	\frac{\Gamma}{\eta^2} \bigg[\int_{\om_0} \int_{ \mathbb{T}^2} G(\vec{x}, \vec{y} )\d  \vec{x} \d \vec{y}
	-
	\int_{\om_0} \int_{ 
		\Omega_0} G(\vec{x}, \vec{y} )\d  \vec{x} \d \vec{y}
	\bigg]  
		\end{align*}
		where
		\begin{align*}
		\frac{\Gamma}{\eta^2} \int_{\om_0} \int_{ \mathbb{T}^2} G(\vec{x}, \vec{y} )\d  \vec{x} \d \vec{y} =O(1),
		\qquad
		\frac{\Gamma}{\eta^2} 
		\int_{\om_0} \int_{ 
		\Omega_0} G(\vec{x}, \vec{y} )\d  \vec{x} \d \vec{y}
			\le O( \eta^2|\log \eta| ). 
		\end{align*}
		%
Therefore,
\begin{align*}
\frac{\Gamma}{\eta^2} \int_{\om_0}  \int_{\mathbb{T}^2} G(\vec{x}, \vec{y} )\d  \vec{x} \d \vec{y} =
\frac{\Gamma}{\eta^2}\bigg[\int_{\mathbb{T}^2} G(\vec{x}, \vec{x}_0)\d \vec{x}+o\bigg(\int_{\mathbb{T}^2} G(\vec{x}, \vec{x}_0)\d \vec{x} \bigg)\bigg] |\Omega_0| .
\end{align*}
%
This follows because we can choose an arbitrary $x_0\in \Omega_0$, and in view of \eqref{bounded perimeter},
as well as the fact that the diameter is bounded from above by the perimeter, 
we have
$\text{diam}\, \Omega_0 $ becomes vanishingly small as $\eta\to 0$.
We then have, for any $\overset{\to}{x}\in \mathbb{T}^2,\overset{\to}{y}\in \Omega_0$,
\begin{align*}
|\overset{\to}{x}-\overset{\to}{x_0}|-\text{diam}\, \Omega_0
\le
|\overset{\to}{x}-\overset{\to}{y}| & \le
|\overset{\to}{x}-\overset{\to}{x_0}|+\text{diam}\, \Omega_0 .
\end{align*}
Then, recalling that
\[G(\overset{\to}{x},\overset{\to}{y})
=-\frac{1}{2\pi}\ln|\overset{\to}{x}-\overset{\to}{y}|+ R (|\overset{\to}{x}-\overset{\to}{y}|)\]
where $R$ is the regular part of $G$ such that $-\Delta R =1$,
we get
\begin{align}
G(\overset{\to}{x},\overset{\to}{y}) - G(\overset{\to}{x},\overset{\to}{x_0})
&=
-\frac{1}{2\pi}  \Big[\ln|\overset{\to}{x}-\overset{\to}{y}| -\ln|\overset{\to}{x}-\overset{\to}{x_0}| \Big]+ R (|\overset{\to}{x}-\overset{\to}{y}|)
-R (|\overset{\to}{x}-\overset{\to}{x_0}|)\notag\\
&=
-\frac{1}{2\pi}  \ln \frac{|\overset{\to}{x}-\overset{\to}{y}|}{|\overset{\to}{x}-\overset{\to}{x_0}|}
+R'({z}) (|\overset{\to}{x}-\overset{\to}{y}|-
|\overset{\to}{x}-\overset{\to}{x_0}|)=:\omega_\eta(\overset{\to}{x},\overset{\to}{y}),\label{how distance changes}
\end{align}
for some $z$ between $|\overset{\to}{x}-\overset{\to}{y}|$ and $|\overset{\to}{x}-\overset{\to}{x_0}|$.
Thus
\begin{align*}
\int_{\Omega_0}  \int_{\mathbb{T}^2}G(\overset{\to}{x},\overset{\to}{y}) \d\overset{\to}{x}\d\overset{\to}{y}
&=
\int_{\Omega_0}  \int_{\mathbb{T}^2}G(\overset{\to}{x},\overset{\to}{x_0}) \d\overset{\to}{x}\d\overset{\to}{y}
+
\int_{\Omega_0}  \int_{\mathbb{T}^2}\omega_\eta(\overset{\to}{x},\overset{\to}{y}) \d\overset{\to}{x}\d\overset{\to}{y}\\
&=
|\Omega_0|  \int_{\mathbb{T}^2}G(\overset{\to}{x},\overset{\to}{x_0}) \d\overset{\to}{x}
+
\int_{\Omega_0}  \int_{\mathbb{T}^2}\omega_\eta(\overset{\to}{x},\overset{\to}{y}) \d\overset{\to}{x}\d\overset{\to}{y},
\end{align*}
where $\omega_\eta(\overset{\to}{x},\overset{\to}{y})\to 0$ a.e., and
is dominated by the integrable function $ |G(\overset{\to}{x},\overset{\to}{y}) |+| G(\overset{\to}{x},\overset{\to}{x_0})|$, and hence
\[\frac{1}{|\Omega_0|}\int_{\Omega_0}  \int_{\mathbb{T}^2}\omega_\eta(\overset{\to}{x},\overset{\to}{y}) \d\overset{\to}{x}\d\overset{\to}{y} \overset{\eta\to 0}{\to} 0,\]
i.e.
\[\frac{1}{|\Omega_0|}\int_{\Omega_0}  \int_{\mathbb{T}^2}\omega_\eta(\overset{\to}{x},\overset{\to}{y}) \d\overset{\to}{x}\d\overset{\to}{y}
=o\bigg( \int_{\mathbb{T}^2}G(\overset{\to}{x},\overset{\to}{x_0}) \d\overset{\to}{x} \bigg) \qquad \text{as }\eta\to 0.\]
%
We thus define a third constant
\begin{align}
\label{c3 value}
c_3=c_3(\Gamma)=\Gamma\int_{\mathbb{T}^2} G(\vec{x}, \vec{x}_0)\d \vec{x},
\end{align}
that will be determined in the $\Gamma$-limit
			 and replace ${e}(m)$ in \eqref{limit functions e}
			 by 
			 ${e_0}(m)={e}(m)+c_3 m$ which will allow us to use the results in \cite{bi1}.
		
		\medskip
		
		%
		Indeed,
		the analogue of $E^{2d}_0$ from \cite{bi1}, should be 
	 of the form
		\begin{align}
		E_0(v):=
		\begin{cases}
				\sum_k  \overline{e_0}(m_k) + 
				\sum_j 
				\tilde{e}_0(m_j)
				& \text{if } v=\sum_k m_k \delta_{\vec{x}_k}+
				\sum_j m_j \delta_{\vec{x}_j},\ m_k,m_j \ge 0,  \text{ with } \{\vec{x}_k,\vec{x}_j \} \text{ all distinct}, 
				\\
				&x_k \text{ on the boundary between the two majority phases},
				\\
				& x_j \text{ not on the boundary between the two majority phases},\\
				+\infty & \text{otherwise},
		\end{cases}
		\label{1st order Gamma limit energy}
		\end{align}
				where
		\begin{align}
		\overline{e_0}(m) &:=\inf \Big\{ \sum_k  e_0  (m_k):m_k\ge 0 ,\sum_k m_k =m \Big\},\qquad
		e_0(m)=c_1 \sqrt{m} + c_2m^2 +c_3m, \label{e0c1}\\
		\tilde{e}_0(m)&:=
		\inf \Big\{ \sum_k [		2\sqrt{\pi} \sqrt{m_k} + c_2m_k^2 +c_3m_k] :m_k\ge 0 ,\sum_k m_k =m \Big\}.
		\end{align}
		The constants $c_1,c_2,c_3$ all depend on the optimal shape of 
		the minority phase,
		which will be described below,
		 with $c_2,c_3$ also depending on $\Gamma$.
		We refer to \eqref{c1 value}, \eqref{c2 value}, \eqref{c3 value} for their specific values.


 	\subsection{Optimal lens}\label{optimal lenses}

		The key difference between our work and \cite{bi1} is that we have two majority phases, and hence the boundary between
		$\Omega_1$ and $\Omega_2$ contributes to the overall energy. This
		was not the case in \cite{bi1}, where only one type of background was present,
		and only the minority phase contributed to the overall energy.
		Therefore, as opposed to \cite{bi1}, we do not expect $\Omega_0$ to be a union of balls.
				
		We illustrate our argument with an example: assume 
the boundary between
the two majority phases has a straight segment $\ell $, then we can place a component
$\Omega_{0,\eta,k}$
to overlap $\ell$. This would add $\H^1(\partial \Omega_{0,\eta,k})
$ to the overall energy, but it also {\em removes} the term $\H^1(\ell\cap\Omega_{0,\eta,k})$.
Thus when placing the minority phase, we remove 
the contribution to the overall perimeter of the part of the boundary 
	between
the two majority phases overlapping with $\Omega_0 $.

        In consequence, it is energetically advantageous for the components $\Omega_{0,\eta,k}$ to overlap with 
	the boundary 
	between
the two majority phases. Further if the straight segment is of length $\gg \eta$,
	then any triple junction should be composed
	of three angles of amplitude $120^\circ$ since all perimeters are weighted equally.
	This suggests that $\Omega_0$ should resemble a lens rather than a ball.
	
	\medskip
	
We note that the assumption that the boundary between
		the two majority phases has a straight segment $\ell $ is not restrictive; indeed
	by \cite[Proposition~2.1]{st}, we know that any local minimum of $E_L$ has 
	$C^{3,\alpha}$ regular
	reduced boundary, for some $\alpha\in (0,1)$, and that these local minima themselves are $C^{1,\alpha}$ regular too.
Thus, around any point lying on the boundary between the two majority phases,  
in the blow up limit, the boundary becomes straight.

	\medskip
	
	 The next result gives the optimal shape that,
	  for a given mass,
	 minimizes the difference between perimeter and diameter:
	 
	 \begin{lemma}
	 	\label{straight line with lens}
	 	\cite{alama2023lens}
	 	Let $\ell$ be a straight line segment of length one. Then the optimal way to place a set $X$ of area $\eta\ll1$
	 	such as to minimize the quantity
	 	\[\mathcal{H}^1(\ell\setminus X) + \mathcal{H}^1(\partial X)\]
	 	is by having $X$ be a lens described in \eqref{optimal lens - shape} below, placed in such a way that $X$ is symmetrical with respect to $\ell$, and the center of $X$ lies on $\ell$.
	 \end{lemma}
	 
	 We will give a proof of this lemma for the sake of completeness and for our specific setting.
	 
	 	\medskip
		
	 First, we show that for any such optimal $X$, the intersections $\ell\cap \partial
	 X$ are triple junctions made of three angles of $120^\circ$.
	 
	 \begin{lemma}
	 	\cite{alama2023lens}
	 	Any such optimal $X$ from Lemma \ref{straight line with lens} is such that the intersections $\ell\cap \partial
	 	X$ are triple junctions made of three angles of $120^\circ$.
	 \end{lemma}
 
 \begin{proof}
 	Since we want to minimize $\mathcal{H}^1(\ell\setminus X) + \mathcal{H}^1(\partial X)$, it is clearly convenient to have
 	$X$ intersecting $\ell$. Then, observe that such $X$ should be convex: its convex hull
 	$\text{conv}(X)$ always satisfies
 	\[\text{diam} (X) =\text{diam}(\text{conv} (X)),\qquad \H^2(X)\le \H^2(\text{conv}(X)),\qquad
 	\H^1(\partial \text{conv}(X))\le \H^1(\partial X), \]
 	with the equality holding only when $X$ is convex.
 	
 	Now that we can work with $X$ convex, the interaction $\ell \cap \partial X$ has at most two points.
 	Denote by $p$ such a junction, and assume that the three angles formed by $\ell$ and $\partial X$ have amplitudes
 	$\alpha,\beta,2\pi-\alpha-\beta$, with at least one of $\alpha$ or $\beta$ not being $120^\circ$ degree.
 	We construct a competitor $X_\varepsilon$ in the following way:
 	denote by $\tau_1$ and $\tau_2$ the tangent lines to $\partial X$, starting from $p$. Choose a ball $B(p,\varepsilon)$
 	of center $p$ and radius $\varepsilon$, and let
 	\[p_0:=\ell \cap  \partial (B(p,\varepsilon) \setminus X),\qquad  p_i:=\tau_i\cap  \partial B(p,\varepsilon) ,\qquad i=1,2.\]
 	Then define
 	\[X_\varepsilon :=(X\setminus B(p,\varepsilon) )\cup \Sigma(p_0,p_1,p_2),\]
 	where $\Sigma(p_0,p_1,p_2)$ denotes the Steiner graph (i.e. the shortest connected tree) connecting $p_0,p_1,p_2$.
 	Since at least one between $\alpha$ and $\beta$ is not $120^\circ$, then
\begin{align}
\label{replace with Steiner graph}
  \H^1( (\ell\cup \partial X)\cap B(p,\varepsilon)  ) - \mathcal{H}^1(\Sigma(p_0,p_1,p_2)) \ge c(\alpha,\beta)\varepsilon, 
\end{align} 
where $c(\alpha,\beta)$ is a constant depending only on $\alpha,\beta$, and strictly positive.
Consequently,
\begin{align}
\label{difference in length}
\mathcal{H}^1(\ell \setminus X )+\mathcal{H}^1(\partial X )=
\mathcal{H}^1(\ell \setminus X_\varepsilon )+\mathcal{H}^1(\partial X_\varepsilon )+ c(\alpha,\beta)\varepsilon.
\end{align}
However, note that $\mathcal{H}^2(X_\varepsilon)$ is not necessarily equal to $\eta$. To overcome this issue,
we observe that the construction in \eqref{replace with Steiner graph} affects only the region inside
$B(p,\varepsilon) $, with area $\pi \varepsilon^2$, thus
$|\mathcal{H}^2(X)-\mathcal{H}^2(X_\varepsilon)|=c'\varepsilon^2$, for some $c'\le \pi$. Then it suffices to apply
the scaling of ratio $\sqrt{1+c'\varepsilon^2} = 1+O(\varepsilon^2)$
 to the whole configuration $\ell \cup X_\varepsilon$: the scaled version
$\sqrt{1+c'\varepsilon^2} X_\varepsilon$ has the required area $\eta$. All the lengths are scaled 
by a factor $\sqrt{1+c'\varepsilon^2}= 1+O(\varepsilon^2)$.
The
set
$\sqrt{1+c'\varepsilon^2} \ell$ might have length other than one:  
if this is the case, we just need to add or remove the length difference, which in any case would not exceed
$\sqrt{1+c'\varepsilon^2}= 1+O(\varepsilon^2)$. Thus, in view of
\eqref{difference in length},
 the above construction decreases 
the overall length. Since we can apply such argument whenever \eqref{difference in length} holds, i.e.
whenever $\alpha$ or $\beta $ is not equal to $120^\circ$, the proof is complete.
 \end{proof}

	 Now we can look for the shape of the optimal lens 
	 (see \cite{alama2023lens2} in the classical setting):

	 \begin{lemma}\label{optimal lens - lemma}
	 	\cite{alama2023lens}
	 	The minimization problem
	 	\begin{align}
	 	\inf &\bigg\{ \int_0^L [\sqrt{1+|f'|^2} + \sqrt{1+|g'|^2}]\d x-L:
	 	 \int_0^L (f-g)
	 	\d x =1,\ f(0)=f(L)=g(0)=g(L)=0,  f, g \in C^1  \bigg\} 
	 	\label{minimization lens}
	 	\end{align}
	 	has a unique solution given by
	 	\begin{align}
	 	\label{optimal lens - shape}
	 	f(x)=-g(x)=\sqrt{3}L \bigg[\sqrt{1-\frac{3}{4}\Big(1-\frac{2x}{L}\Big)^2}-\frac{1}{2}\bigg],
	 	\end{align}
	 	with $L$ satisfying
	 	\begin{align}
	 	\label{mass constraint on L}
	 	2\sqrt{3}L \int_0^L\bigg[\sqrt{1-\frac{3}{4}\Big(1-\frac{2x}{L}\Big)^2}-\frac{1}{2}\bigg]=1.
	 	\end{align}
	 \end{lemma}
	 
	 Again we include a simple proof that applies to our setting, which assumes connectedness.
%
	 \begin{proof}
	 	Symmetry considerations immediately give $f=-g$. 
	 			Any solutions $f,g$ of the minimization problem \eqref{minimization lens} also satisfy
	 	\begin{align}
	 	\label{120 degree angles - formula}
	 	f'(0)=-f'(L)=-g'(0)=g'(L)=\sqrt{3}. 
	 	\end{align}
	 	Thus \eqref{minimization lens} can be reduced to
	 	\begin{align}
	 \inf &\bigg\{ 2\int_0^L \sqrt{1+|f'|^2} \d x-L:
	 2\int_0^L f
	 \d x =1, f(0)=f(L)=0,f'(0)=-f'(L)=\sqrt{3} \bigg\} 
	 \label{minimization lens - only f}
	 \end{align}
	 	Let $\vph\in C^\8([0,L];\mathbb{R})$ be a test function such that		
	 	\[ \vph(0)=\vph(L)=\vph'(0)=\vph'(L)=0,\qquad \int_0^L\vph \d x=0,\]
	 	and direct computations give
	 	\begin{align*}
	 	\lim_{\vep\to 0}& \frac{1}{\vep}\int_0^L [\sqrt{1+|f'+\vep \vph'|^2} -\sqrt{1+|f'|^2} ]\d x=
	 	\lim_{\vep\to 0}
	 	\int_0^L \frac{2 f'\vph' + \vep 
	 		|\vph'|^2}{\sqrt{1+|f'+\vep \vph'|^2} +\sqrt{1+|f'|^2}}
	 	\d x\\
	 	&=\int_0^L \frac{ f'\vph' }{\sqrt{1+|f'|^2}}
	 	\d x = -\int_0^L \bigg[\frac{ f' }{\sqrt{1+|f'|^2}}\bigg]' \vph\d x.
	 	\end{align*}
	 	Since $\vph$ has to satisfy $\int_0^L\vph \d x=0$, the criticality 
	 	condition on $f$ reads
	 	\begin{align*}
	 	\bigg[\frac{ f' }{\sqrt{1+|f'|^2}}\bigg]' = A \in \R \Lra 
	 	\frac{ f' }{\sqrt{1+|f'|^2}} = Ax+B, \qquad A,B\in \R.
	 	\end{align*}
	 	Conditions $f'(0)=-f'(L)=\sqrt{3}$ give
	 	\begin{align*}
	 	B=	\frac{ f'(0) }{\sqrt{1+|f'(0)|^2}} =\frac{\sqrt{3}}{2},\qquad
	 	AL+B = \frac{ f'(L) }{\sqrt{1+|f'(L)|^2}} =-\frac{\sqrt{3}}{2}
	 	\Lra A= -\frac{\sqrt{3}}{L}.
	 	\end{align*}
	 	Hence 
	 	\begin{align*}
	 	\frac{ f' }{\sqrt{1+|f'|^2}} = \frac{\sqrt{3}}{2}\Big(1-\frac{2x}{L}\Big)
	 	\Lra f'(x)^2=\frac{3}{4}\Big(1-\frac{2x}{L}\Big)^2(1+f'(x)^2)
	 	\Lra f'(x)^2 =\frac{\frac{3}{4}(1-\frac{2x}{L})^2}{1-\frac{3}{4}(1-\frac{2x}{L})^2},
	 	\end{align*}
	 	thus
	 			\[f'(x) =\frac{\sqrt{3}}{2}\frac{1-\frac{2x}{L}}{ \sqrt{1-\frac{3}{4}(1-\frac{2x}{L})^2}}.\]
	 			Using $f(0)=f(L)=0$, we get
	 			\[ f(x) =\sqrt{3}L \bigg[\sqrt{1-\frac{3}{4}\Big(1-\frac{2x}{L}\Big)^2}-\frac{1}{2}\bigg].\]
	 			The mass constraint condition then forces $L$ to satisfy \eqref{mass constraint on L}.
	 \end{proof}
	 
The set described in Lemma \ref{optimal lens - lemma} is the shape optimum for Lemma \ref{straight line with lens},
and will be denoted by $\Lambda^*$. Any scaled copy of it, i.e. sets of the form $r\Lambda^*$, $r>0$,
will be referred to as {\em optimal lens}. This also gives us an indication about the value of $c_1$ from \eqref{limit functions e}: assuming that, under blow up, the boundary between the two majority phases is straight,
hence it can be approximated by a tangent line $\ell$, then the optimal way is to add an optimal lens placed to be
symmetrical with respect to $\ell$. This adds 
\[c(m) 2\int_0^L \sqrt{1+|f'|^2} \d x\]
to the perimeter,
with $f$ (resp. $L$) defined in \eqref{optimal lens - shape} (resp. \eqref{mass constraint on L}), while removing $c(m) L$,
for some $c(m)$ that depends on the mass $m$ of the optimal lens.
		Here, in 2D, $c(m)=\sqrt{m}$.
	Such scaling factor $c(m)$ is required since the optimal lens defined in Lemma
	\ref{optimal lens - lemma} has unit mass.
Therefore, we should have
\begin{align}
\label{c1 value}
c_1=2\int_0^L \sqrt{1+|f'|^2} \d x - L.
\end{align}
The geometry of optimal lenses also give us a guess about the value of 
 $c_2=c_2(\Gamma)$. It comes from
 the self-interaction of the minority phase, i.e., the term:
\[\frac{1}{\eta^3 |\log \eta|}\int_{\Omega_0 } \int_{\Omega_0
} G(\vec{x}, \vec{y} )\d \vec{x} \d \vec{y},\]
which, using $\Omega_{0,\eta}=\bigsqcupdot_{k} \Omega_{0,\eta ,k} $, can be split into 
\begin{align*}
\frac{1}{\eta^3 |\log \eta|} \bigg[\sum_k \int_{\Omega_{0,\eta ,k} } \int_{ \Omega_{0,\eta ,k}
} G(\vec{x}, \vec{y} )\d \vec{x} \d \vec{y} +
\sum_{k\neq l} \int_{\Omega_{0,\eta ,k}} \int_{ \Omega_{0,\eta ,l}
} G(\vec{x}, \vec{y} )\d \vec{x} \d \vec{y}
\bigg] 
\end{align*}
Since, as $\eta\to 0$, the mutual distance between different $\Omega_{0,\eta ,k}$ is of order $O(1)$, while the diameter of each single 
$\Omega_{0,\eta ,k}$ is of order $O(\eta)$, we get that
\[\sum_{k\neq l} \int_{\Omega_{0,\eta ,k}} \int_{ \Omega_{0,\eta ,l}
} G(\vec{x}, \vec{y} )\d \vec{x} \d \vec{y}
= o\bigg(  \sum_k \int_{\Omega_{0,\eta ,k}} \int_{ \Omega_{0,\eta ,k}
} G(\vec{x}, \vec{y} )\d \vec{x} \d \vec{y}
 \bigg).\]
Thus the leading order are given by the self interactions, and hence
\begin{align}
\label{c2 value}
c_2:= \begin{cases}
\Gamma_{00} \, c_{2,B} & \text{if the corresponding mass is inside one background type},\\
\Gamma_{00} \, c_{2,L} & \text{if the corresponding mass is on the boundary between the two backgrounds},
\end{cases}
\end{align}
where
\begin{align*}
c_{2,B}&:=\int_{B_1\times B_1
} G(\vec{x}, \vec{y} )\d \vec{x} \d \vec{y},\qquad B_1:=\text{ball of unit mass},\\ 
c_{2,L}&:=\int_{\Lambda_1\times \Lambda_1
} G(\vec{x}, \vec{y} )\d \vec{x} \d \vec{y},\qquad \Lambda_1:=\text{optimal lens of unit mass}.
\end{align*}
It is also known that the ball is the worst case for the interaction term \cite{abcase}, hence
$c_{2,B}>c_{2,L}$. This, combined with the fact that $c_1<2\pi$,
indicate that optimal lenses are energetically preferable to balls, and hence,
in optimal configurations,
all type 0 masses tend to lie on the boundary between the two backgrounds.

		
 \subsection{Proof of the Theorem \ref{1st order Gamma convergence}} \label{proof of thm 2.2}		
		
	 The main result of this section is:
\begin{theorem}
	\label{1st order Gamma convergence}
	It holds $\frac{E_\eta-E_L}{\eta} \overset{\Gamma}{\to} E_0 $.
\end{theorem}

	 Again, the main difference with \cite{bi1} is that we have two different
	 backgrounds, thus a limit configuration will have a partition
	 \begin{align}
	 \label{limit configuration - background}
	 \bigsqcupdot_{k} \Omega_{1,k} \sqcupdot 	 \bigsqcupdot_{k} \Omega_{2,k} 
	 ,\qquad \text{all } \Omega_{i,k} \text{ mutually disjoint},
	 \end{align}
	 where the $\Omega_{i,k}$ denote regions of type $i$ constituent ($i=1,2$),
	 while the minority phase is represented by
	 \begin{align}
	 \label{limit configuration - dirac masses}
	 \sum_k m_k \delta_{\vec{x}_k},\qquad \text{all } \vec{x}_k \text{ disjoint}.
	 \end{align}
	 This introduces a slight complication: although the energy
	 $\frac{E_\eta-E_L}{\eta}$ does not see $\Omega_{i,k}$, $i=1,2$,
	 hence the partition in \eqref{limit configuration - background} is not influential
	 to $\frac{E_\eta-E_L}{\eta}$,
	 the presence of a boundary between type I and II constituents, combined with the fact that,
	 a priori, it is unclear if any 
	 $\Omega_{i,k}$ has interior part, makes the construction
	 of the recovery sequence in the $\Gamma-\limsup$ part geometrically
	    more delicate.
	 
	 \medskip
	 
	  The compactness part 
	 is argued like in \cite{bi1}
	 and follows from a concentration-compactness result, while the $\Gamma-\liminf$ part is 
	 similar to \cite{bi1} Theorem 6.1.
	 The $\Gamma-\limsup$ part is achieved by substituting Dirac masses with 
	 balls if they lie inside one type of background, and with
	 optimal lenses if they lie on the boundary between types I and II backgrounds.



 \begin{proof}
	 	The proof is split into three parts.
	 	
	 	\medskip
	 	
	 	{\em Compactness.} 
		Note that the extra term in the energy (compared to the energy in \cite{bi1}) is
			the interaction of the minority phase with the two backgrounds, i.e., there is a new term
	 	\[ \frac{\Gamma_{i0} }{\eta}\int_{\om_i} \int_{ \om_0} G(\vec{x}, \vec{y} )\d \vec{x} \d \vec{y}, \]
	 	 which
	 	 is bounded from below, hence
	 	compactness follows from the same arguments in \cite[Lemma~5.1]{bi1}.
	 	
	 	\medskip
	 	
	 	$\Gamma-\liminf$ {\em inequality.}
	 	Note that 
	 	the extra term (compared to \cite{bi1}) 
		\[  \frac{\Gamma_{i0} }{\eta}\int_{\om_i} \int_{ \om_0} G(\vec{x}, \vec{y} )\d \vec{x} \d \vec{y} \]
	 	is linear in the minority phase.
	 	Therefore, the proof follows the same computations from the proof of the
	 		$\Gamma-\liminf$ inequality part of \cite[Theorem~6.1]{bi1}.
	 		
	 	\medskip
	 	
	 	$\Gamma-\limsup$ {\em inequality.} Consider an arbitrary limit configuration, with background
	 	given by a disjoint union
	 	 \begin{align}	 \label{limit configuration - background Gamma limsup}
	 		\bigsqcupdot_{k} \Omega_{1,k} \sqcupdot 	 \bigsqcupdot_{k} \Omega_{2,k} 
	 		,\qquad \text{all } \Omega_{i,k} \text{ mutually disjoint},
	 	\end{align}
	 	and type 0 constituent given by $v=\sum_k m_k \delta_{\vec{x}_k}$, where
	 	$\sum_k m_k=M$, and
	 	 $\vec{x}_k$ are distinct points in $\mathbb{T}^2$.
	 	We construct a recovery sequence in the following way:
	 	\begin{enumerate}
	 			 		
	 		\item If the sum $v=\sum_k m_k \delta_{\vec{x}_k}$ is infinite, then we can make it finite
	 		at the cost of an arbitrarily small error
			by e.g. removing the infinite tail whose total mass (and perimeter) does not exceed an arbitrarily chosen small value.
			Similarly, if the unions in 
	 		\eqref{limit configuration - background Gamma limsup} are infinite, we can make them finite at a cost of an arbitrarily small error again by removing the infinite tail.
	 		 Thus we assume that both unions in \eqref{limit configuration - background Gamma limsup}, and the sum
	 		$v= \sum_{k=1}^{N}  m_k \delta_{\vec{x}_k}$ are finite. In this way, the minimum distance
	 		between two such masses is positive, i.e. 
	 		\[ \inf_{i\neq j} |\vec{x}_i - \vec{x}_j|>0 \]
			and we can use the results of section 3.1.

	 		\item Without loss of generality, let $\vec{x}_1,\cdots, \vec{x}_h$ be 
	 		in the interior of one type of background, while 
	 		$\vec{x}_{j}$, $ N \geq j \ge h+1$, lie on the boundary between the two background types.
	 		Use a ball $B_k:=B(\vec{x}_k, \eta \sqrt{m_k/\pi})$ to approximate each 
	 		$m_k \delta_{\vec{x}_k}$, $k=1,\cdots,h$. That is, define 
	 		\begin{align}
	 		\label{v eta - ball part}
	 			 		v_{\eta, B}:=\frac{1}{\eta^2}\sum_{k=1}^{h} \rchi_{B(\vec{x}_k, \sqrt{m_k/\pi})} .
	 		\end{align}

	 		\item We cannot yet use optimal lenses to approximate $\vec{x}_{j}$, $j\ge h+1$, since the construction of
	 		optimal lenses require a locally straight boundary. Therefore, we first approximate the background sets
	 		$\Omega_{j,k}$, $j=1,2$, with smooth sets $\Omega_{j,k}'$, $j=1,2$, such that 
	 		$\|\Omega_{j,k}\|_{L^1(\mathbb{T}^2)}=\|\Omega_{j,k}'\|_{L^1(\mathbb{T}^2)}$ for all $j,k$, and
	 		$\chi_{\Omega_{j,k}'}$ all converge to 
	 		$\chi_{\Omega_{j,k}}$ in the topology of $BV(\mathbb{T}^2)$. This is possible since smooth functions are dense
	 		in $BV(\mathbb{T}^2)$. Then, we can approximate each $\Omega_{j,k}'$, $j=1,2$, with another
	 		sequence $\Omega_{j,k}'$, $j=1,2$, 
	 		 such that 
	 		 		$\|\Omega_{j,k}''\|_{L^1(\mathbb{T}^2)}=\|\Omega_{j,k}'\|_{L^1(\mathbb{T}^2)}$ for all $j,k$, and
	 		 the following condition is satisfied: around each $\vec{x}_j$, $j\ge h+1$, the boundary
	 		  contains a straight line
	 		 segment of length $\sqrt{\eta}$ centered around $\vec{x}_j$. This is again an vanishing perturbation
	 		 as $\eta \to 0$, thus it does not affect the $\Gamma-\limsup$ inequality.
	 		
	 		\item Now we can
	 			 		place an optimal lens
	 		$X_j$ of mass $\eta^2 m_j$ centered around $\vec{x}_j$, and with diameter lying
	 		on the boundary between the two background types. 
	 		Note that this reduces the overall length of the boundary between types I and II backgrounds, since
	 		part of it is now covered by the $X_j$.
	 		This construction is possible (for all sufficiently small $\eta$) since $X_j$ will have diameter $O(\eta)$, while in the previous step we ensured that the boundary contains a straight line
	 		segment of length $\sqrt{\eta}$ centered around $\vec{x}_j$. Then define 
	 		\[v_{\eta,L} := \frac{1}{\eta^2} 
	 		\sum_{j\ge h+1} \rchi_{X_j},\]
	 		and let our approximating sequence be
	 		$v_\eta:=v_{\eta,B}+v_{\eta,L}$.
	 	\end{enumerate}
 	Since the constant $c_1$ from the definition of $\overline{ e_0} $ in (\ref{e0c1}) is less than $2\pi$,
	the same arguments from the proof of
	Theorem~6.1 from \cite{bi1} gives the $\Gamma-\limsup$ inequality.
 		 \end{proof}
	
	{\bf Remark.}
	From the proof of the $\Gamma-\limsup$ inequality we can see that, since 
	the limit energy $E_0$ is invariant with respect to positions of $\vec{x}_k$, and
	$c_1<2\pi$, $c_{2,L}<c_{2,B}$, having even a single $\vec{x}_k$
	not on the boundary between the two background types makes the
	$\Gamma-\limsup$ inequality {\em strict}, i.e. the limit configuration cannot be a global minimum. 
		

\vspace{0.2cm}

\section{Uniform Masses} \label{uniform}

In this section we are going to prove that, in optimal configurations, the type 0 constituent is spread over a finite number of equally sized masses.
The main result is:
\begin{theorem}
	\label{uniform lenses theorem}
	Let $M, \Gamma_{ij}$ be given. Then, in any optimal configuration, the type 0 constituent is
		spread over at most
	$1+ M(\frac{c_1}{8c_2})^{-2/3}$
	identically sized masses.
\end{theorem}

\begin{proof}
	The proof is split over several steps. 
		Let an arbitrary optimal configuration be given, with $v=\sum_k m_k \delta_{x_k}$, $\sum_k m_k=M$,
		describing the type 0 constituents. If there is only one such $m_k$, then there is nothing to prove.
		Thus we will assume, in the rest of the proof, the existence of at least two such masses $m_k$.
		
		\medskip
		
		{\em Step 1. Finite mass splitting.}
		First, we note that 
		$$e_0(m)=c_1 \sqrt{m} + c_2m^2 +c_3m,
		$$ 
		is concave for small $m$, hence in
		\begin{align*}
		\overline{e_0}(m):=\inf \Big\{ \sum_k e(m_k):m_k\ge 0 ,\sum_k m_k =m \Big\},\qquad
		e_0(m)= \underbrace{ c_1 \sqrt{m} + c_2m^2}_{=:f(m)} +c_3m,
		\end{align*}
		the infimum is never an infinite sum, but always a finite one. That is,
		we always have
		\[\overline{e_0}(m)=\inf \Big\{ \sum_k e(m_k):m_k\ge 0 ,\sum_k m_k =m \Big\}=\sum_{k=1}^N e(m_k)
		=c_3m+\sum_{k=1}^N f(m_k),
		\]
		for some suitable $m_1,\cdots,m_N$.
		Similarly,
		 the limit energy
	can be rewritten as
		\begin{align*}
	E_0(v)=
	\sum_k \overline{e_0}(m_k)=\sum_k\Big[c_3m_k+ \sum_{j=1}^{N(k)} f(m_{k,j}) \Big]=c_3M+
	\sum_k \sum_{j=1}^{N(k)} f(m_{k,j}) ,
	\end{align*}
	for suitable (still potentially infinite) mass splitting $m_{k,j}$, $k\geq 1$, $j=1,\cdots,N(k)$. Using again the concavity
	of $f$ for small $m$, we get that any $v$ minimizing $E_0$ must be of the form
	$v=\sum_{k=1}^{K} m_k \delta_{\vec{x}_k}$, $\sum_k m_k=M$, for some suitable $K$. Since we proved before that any optimal splitting realizing the infimum for $\overline{ e_0}$ must also be finite, we infer that, 
		for any optimal configuration,
	the energy is also a finite sum. Hence, by relabeling the masses, we have
	\[E_0(v)=c_3M+E(v),\qquad E(v):=
	\sum_{k=1}^{N}f(m_{k}),\]
	for some suitable $N$, and mass splitting $m_1,\cdots,m_N$. Clearly, as $M$ is given, minimizing $E_0$
	is equivalent to minimizing $E$. 
	
	\medskip
	
	{\em Step 2. Uniform upper bound on optimal mass splittings.} 
	We note that the $N$ from before can depend on
	the optimal configuration 
	 $v$ itself. To show that $N$ admits an upper bound which depends only on the parameters of the minimization, namely $M$ and
	 $\Gamma_{ij}$, we again rely on the concavity of $f$ for small $m$: its second derivative
	 is $-\frac{c_1}{4}m^{-3/2}+2c_2$, hence it is concave up to $(\frac{c_1}{8c_2})^{2/3}$, and convex beyond this value.
	 As such, any $v$ minimizing $E$ can have at most one mass in $[0,(\frac{c_1}{8c_2})^{2/3}]$. Therefore, an upper bound
	 for the number of masses $N$ is 
	 \[1+ M\Big(\frac{c_1}{8c_2}\Big)^{-2/3},\]
	 which clearly depends only on $M,\Gamma_{ij}$.
	 
	 \medskip
	 
	 {\em Step 3. Two candidate optimal configurations.} Using the fact that $f$ is convex
	 for $m> (\frac{c_1}{8c_2})^{2/3}$, in any optimal configuration all the masses in 
	 $[(\frac{c_1}{8c_2})^{2/3},M]$ must be all equal.
	  Thus there are essentially two classes of
	 candidates for
	 optimal configurations:
	 \begin{enumerate}
	 	\item one mass $m\in (0,(\frac{c_1}{8c_2})^{2/3})$, and the remaining $N-1$ ones equal to $\frac{M-m}{N-1}\in [(\frac{c_1}{8c_2})^{2/3},M]$;
	 	\item all masses being in $[(\frac{c_1}{8c_2})^{2/3},M]$. In this case, all of them are equal to
	 	$M/N$.
	 \end{enumerate}
%
%
We need to show that the first case is never minimizing.
In this case the energy depends on $N$ and $m$, and can be written as
\begin{align*}
E(N,m)& = c_1\sqrt{m} +c_2m^2 + (N-1)\bigg[ c_1 \sqrt{\frac{N-m}{N-1}} + c_2\Big(\frac{N-m}{N-1}\Big)^2 \bigg]\\
&=c_1\sqrt{m} +c_2m^2 + c_1 \sqrt{(N-m)(N-1)}+c_2\frac{(N-m)^2}{N-1}.
\end{align*}
The partial derivative in $m$ is
\begin{align*}
\frac{\partial}{\partial m}E(N,m) &=\frac{c_1}{2}\bigg[\frac{1}{\sqrt{m}}-\sqrt{ \frac{N-1}{N-m} } \bigg]+2c_2 \bigg[
m-\frac{N-m}{N-1}\bigg] .
\end{align*}
Note that $\frac{N-m}{N-1}$ must lie in the convex region of $f$, while
$m$ lies in its concave region. As such, $\frac{N-m}{N-1}\ge m$, which means that
$\frac{\partial}{\partial m}E(N,m)\le 0$, with equality reachable only when 
$m = (\frac{c_1}{8c_2})^{2/3}$. But in this case, all the $N$ masses would belong to
$[(\frac{c_1}{8c_2})^{2/3},M]$, i.e. the convex region of $f$. The following dichotomy arises:
\begin{enumerate}
	\item either all masses are equal to  $(\frac{c_1}{8c_2})^{2/3}$;
	\item if not, we would have one mass equal to $(\frac{c_1}{8c_2})^{2/3}$, and the remaining $N-1$ ones lying in
	$((\frac{c_1}{8c_2})^{2/3},M]$. Again, the convexity of $f$ here
	makes the all uniform configuration energetically preferable.
\end{enumerate}
Thus in both cases having all uniform masses is energetically preferable.
The proof is thus complete.
\end{proof}

\begin{corollary}
	Any optimal configuration, where type 0 constituent is spread over at least two masses, will have either 
	\[ \bigg\lfloor M\Big(\frac{c_1}{2c_2}\Big)^{-2/3} \bigg\rfloor\qquad \text{or} \qquad 
	\bigg\lceil
	M\Big(\frac{c_1}{2c_2}\Big)^{-2/3} \bigg\rceil\]
	type 0 constituent masses of equal size. 
	Here, 
	\begin{align*}
	\lfloor x  \rfloor :=\max \{n\in\mathbb{Z} :n\le x \},\qquad
	\lceil x \rceil:=\min \{n\in\mathbb{Z} :n\ge x \},
	\end{align*}
	are the floor and ceiling functions respectively.
\end{corollary}

We remark that the condition of having at least two masses is not very strong: by comparing the configuration
with only one mass equal to $M$, and the one with two masses equal to $M/2$, their energies are
\[ c_1\sqrt{M}+c_2M^2 \qquad\text{and}\qquad c_1\sqrt{2M}+\frac{c_2}{2}M^2\]
respectively, thus the latter becomes preferable when $M>(2(\sqrt{2}-1)\frac{c_1}{c_2})^{2/3}$.

\begin{proof}
	In Theorem \ref{uniform lenses theorem} we proved that, in optimal configurations, all masses $m_k$ are equal, 
	thus
	we can find the optimal value. 
	
	\medskip
	
	The energy now depends only on $N$, and can be written as
	\begin{align*}
	E(N)= Nf\Big(\frac{M}{N} \Big) = N\bigg[ c_1 \sqrt{\frac{M}{N}}+c_2\Big(\frac{M}{N}\Big)^2\bigg].
	\end{align*}
	Its derivatives are
	\begin{align*}
	E'(N)&= \frac{M( c_1N -2Mc_2 \sqrt{\frac{M}{N}} )}{2N^2\sqrt{\frac{M}{N}}},\qquad
	E''(N)= \frac{M^2(8c_2(\frac{M}{N})^{3/2}-c_1  )}{4(\frac{M}{N})^{3/2}N^3}.
	\end{align*}
	Note the second derivative $E''$ changes sign exactly when
	$8c_2(\frac{M}{N})^{3/2}-c_1 =0$,
	i.e. when the masses are equal to $(\frac{c_1}{8c_2})^{2/3}$. In other words, $E$ and $f$ are both concave
	before $(\frac{c_1}{8c_2})^{2/3}$, and convex afterwards.
	Since all the masses are equal, they are at least $(\frac{c_1}{8c_2})^{2/3}$, i.e. the second derivative
	$E''$ is nonnegative.
	Direct computations give that $E'(N)= 0$
	at 
	\[N=M\Big(\frac{c_1}{2c_2}\Big)^{-2/3},\]
	where each mass is equal to
	\[\frac{M}{N}=\Big(\frac{c_1}{2c_2}\Big)^{2/3},\]
	which lies in the convexity region of $E$, thus corresponding to a minimum for $E$, provided that
	\[N=M\Big(\frac{c_1}{2c_2}\Big)^{-2/3}\in\mathbb{N}.\]
	On the other hand, if
	\[N=M\Big(\frac{c_1}{2c_2}\Big)^{-2/3}\notin\mathbb{N},\]
	then we note that:
	\begin{enumerate}
		\item $\frac{M}{ \lfloor N \rfloor}$ lies in the convexity region of $E$
		since 
		$$\frac{M}{ \lfloor N \rfloor}> \frac{M}{N}=\Big(\frac{c_1}{2c_2}\Big)^{2/3}; $$

		\item $\frac{M}{ \lceil
			N \rceil}$ 	 also lies in the convexity region of $E$
				since $N\ge 2$, $\frac{ \lceil N \rceil}{N}< \frac{3}{2}$, hence
		$$\frac{M}{ \lceil
	 N \rceil}> \frac{2}{3}\frac{M}{N}=\frac{2}{3}\Big(\frac{c_1}{2c_2}\Big)^{2/3} > \Big(\frac{c_1}{8c_2}\Big)^{2/3}.$$
	\end{enumerate}
	Thus the two configurations with $\frac{M}{ \lfloor N \rfloor}$ and $\frac{M}{ \lceil
		N \rceil}$
	identical masses are acceptable candidate optimal configurations, since both
	$\frac{M}{ \lfloor N \rfloor}$ and $\frac{M}{ \lceil
		N \rceil}$
	 lie in the convexity region of $E$. Since we showed that $E'=0$ only at $M(\frac{c_1}{2c_2})^{-2/3}$,
	 i.e. $E$ is strictly decreasing before $M(\frac{c_1}{2c_2})^{-2/3}$
	 and increasing afterwards,
	the optimal configuration can be made only by either $\frac{M}{ \lfloor N \rfloor}$ or $\frac{M}{ \lceil
		N \rceil}$ identical masses.

\end{proof}


\section{Negative interaction coefficients}

Throughout this section, we allow all interaction coefficients to be negative,
and show zeroth- and first-order $\Gamma$-convergence results when such coefficients are sufficiently
small.

\begin{lemma}\label{perimeter controls interaction}
	Given two (sufficiently regular) sets $E_1, E_2\subseteq \mathbb{T}^2$,
	for which $\mathbf{1}_{E_1}, \mathbf{1}_{E_2}\in BV(\mathbb{T}^2)$, and
	 satisfying 
	$ |E_2|\le  |E_1|$,
	 it holds
	 \begin{align}
	 \frac{1}{2}\sum_{i=1}^{2}\Peri(E_i) &\ge  \frac{c_{isoper, \mathbb{T}^2}}{2}  \sum_{i=1}^{2}|E_{i}|^{1/2},
	 \label{perimeter controls interaction - perimeter estimate}\\
	 	-|E_2||E_1|\sup_{\vec{x}, \vec{y}\in \mathbb{T}^2} |R(\vec{x}, \vec{y})|\le
	 \int_{E_1 \times E_2 } G_{}(\vec{x}, \vec{y})  d\vec{x} d\vec{y} & \le \bigg(\frac{1-\ln |E_1|+\ln \pi}{4\pi}  +\sup_{\vec{x}, \vec{y}\in \mathbb{T}^2} |R(\vec{x}, \vec{y})|\bigg)
	 |E_2| |E_1|	, \label{perimeter controls interaction - interaction estimate}
	 \end{align}
	 where $c_{isoper, \mathbb{T}^2}>0$ is the isoperimetric constant of $\mathbb{T}^2$, and
	 $R$ denotes the regular part of Green's function $G$.
	 
	 As a consequence, since 
	 	$|E_1|\le |\mathbb{T}^2|\le 1$, we get
	 	$\frac{1-\ln |E_1|+\ln \pi}{{4\pi}}\ge 0$, which in turn gives
	 	\begin{align}
	 		\bigg|\int_{E_1 \times E_2 } G_{}(\vec{x}, \vec{y})  d\vec{x} d\vec{y}\bigg| &
	 		\le  \bigg(\frac{1-\ln |E_1|+\ln \pi}{{4\pi}}  +\sup_{\vec{x}, \vec{y}\in \mathbb{T}^2} |R(\vec{x}, \vec{y})|\bigg)
	 		|E_2| |E_1|\label{perimeter controls interaction - interaction estimate - with absolute values}
	 	\end{align}
\end{lemma}

\begin{proof}
	Assertion \eqref{perimeter controls interaction - perimeter estimate} follows directly from the isoperimetric inequality.
	The interaction term satisfies
	\begin{align*}
	\int_{E_1 \times E_2 } G_{}(\vec{x}, \vec{y})  d\vec{x} d\vec{y} 
	&=
	\int_{E_1 \times E_2  } -\frac{1}{2\pi}\ln |\vec{x}-\vec{y}|  d\vec{x} d\vec{y} 
	+\int_{E_1 \times E_2  } R(\vec{x},\vec{y} ) d\vec{x} d\vec{y} 
	\\
	&\le
	|E_2|\int_{ B(0,\sqrt{ |E_1|/\pi})}  -\frac{1}{2\pi}\ln |\vec{x}| \d \vec{x}
	+|E_1||E_2|\sup_{\vec{x}, \vec{y}\in \mathbb{T}^2} |R(\vec{x}, \vec{y})|\\
	&=
	|E_2|\int_0^{ \sqrt{ |E_1|/\pi} }  -r\ln r \d r
	+|E_2||E_1|\sup_{\vec{x}, \vec{y}\in \mathbb{T}^2} |R(\vec{x}, \vec{y})|\\
	&= \frac{1}{2} |E_2| 
	\bigg[ \frac{r^2}{2} -r^2\ln r \bigg] \bigg|_{r=0}^{r= \sqrt{ |E_1|/\pi}} 
	+|E_2||E_1|\sup_{\vec{x}, \vec{y}\in \mathbb{T}^2} |R(\vec{x}, \vec{y})|
	\\
	&=
	|E_2| |E_1|
	\frac{1-\ln |E_1|+\ln \pi}{ 4\pi }  
	+|E_2||E_1|\sup_{\vec{x}, \vec{y}\in \mathbb{T}^2} |R(\vec{x}, \vec{y})|,
		\end{align*}
	and the upper bound in \eqref{perimeter controls interaction - interaction estimate} is proven.
	 For the lower bound, it suffices to notice that
	\begin{align*}
		\int_{E_1 \times E_2 } G_{}(\vec{x}, \vec{y})  d\vec{x} d\vec{y} 
		&=
		\int_{E_1 \times E_2  } -\frac{1}{2\pi}\ln |\vec{x}-\vec{y}|  d\vec{x} d\vec{y} 
		+\int_{E_1 \times E_2  } R(\vec{x},\vec{y} ) d\vec{x} d\vec{y} 
		\\
		&\ge
		|E_2|\int_{ B(0,\sqrt{ |E_1|/\pi})}  -\frac{1}{2\pi}\ln |\vec{x}| \d \vec{x}
		-|E_1||E_2|\sup_{\vec{x}, \vec{y}\in \mathbb{T}^2} |R(\vec{x}, \vec{y})|\\
		&=
		|E_2|\int_0^{ \sqrt{ |E_1|/\pi} }  -r\ln r \d r
		-|E_2||E_1|\sup_{\vec{x}, \vec{y}\in \mathbb{T}^2} |R(\vec{x}, \vec{y})|\\
		&= \frac{1}{2} |E_2| 
		\bigg[ \frac{r^2}{2} -r^2\ln r \bigg] \bigg|_{r=0}^{r= \sqrt{ |E_1|/\pi}} 
		-|E_2||E_1|\sup_{\vec{x}, \vec{y}\in \mathbb{T}^2} |R(\vec{x}, \vec{y})|
		\\
		&=
		\underbrace{|E_2| |E_1|
		\frac{1-\ln |E_1|+\ln \pi}{ 4\pi }  }_{\ge 0}
		-|E_2||E_1|\sup_{\vec{x}, \vec{y}\in \mathbb{T}^2} |R(\vec{x}, \vec{y})|,
	\end{align*}
	concluding the proof.
\end{proof}

 \begin{theorem}\label{zeroth order gamma limit - negative interaction coefficients}
	The conclusion of Lemma \ref{zeroth order gamma limit} can be extended
	to allow the interaction coefficients to be negative, as long as these satisfy
	\begin{align}
	\max_{i=1,2} (  |\gamma_{i1}|+ |\gamma_{i2}|  ) <
	\frac{  \sqrt{2} c_{isoper, \mathbb{T}^2}}{c_{3,R,\mathbb{T}^2}}  ,\qquad c_{3,R,\mathbb{T}^2}:=\frac{1+\ln 3+\ln \pi}{4\pi}  +\sup_{\vec{x}, \vec{y}\in \mathbb{T}^2}|R(\vec{x}, \vec{y})|.
	\label{zeroth order gamma limit - negative interaction coefficients - conditions}
	\end{align}
\end{theorem}

\begin{proof}

	Consider a sequence $\eta_n\to 0$, and 
	\[ u_{n,i}:\mathbb{T}^2 \longrightarrow \{0,1\},\ 	\sum_{i = 0}^2 u_{n,i} =1,
	\ 	\Omega_{n,i}:=\supp u_{n,i},
	\ |\Omega_{n,0}| = \eta_n^2 M,
	\ \frac{1}{3}\le |\Omega_{n,1}|=|\Omega_{n,2}|=\frac{1-\eta_n^2M}{2},
	\qquad i=0,1,2. \]	
	We notice that 
	since as  $\eta \searrow 0$, the masses $|\Omega_{n,1}|=|\Omega_{n,2}| \nearrow 1/2 $, we can replace $1/3$ with any quantity strictly less than $1/2$.

{\em Compactness.} When the energies are uniformly bounded,
clearly
$u_{n,0} \to 0$ strongly in $L^p(\mathbb{T}^2)$ for all $p<+\infty$, and
we aim to show the
existence of $u_i\in BV(\mathbb{T}^2)$ ($i=1,2$) such that, up to subsequence
which will not relabeled,
$u_{n,i} \to u_i$ weakly in $BV(\mathbb{T}^2)$. The idea is that, when the interaction coefficients
are sufficiently small, they are dominated by the perimeter, 
so the energy still controls the perimeter, which in turn implies
boundedness in $BV(\mathbb{T}^2)$.
First we note that interactions involving minority phases are vanishing, i.e.
\begin{align*}
\int_{\mathbb{T}^2 \times \mathbb{T}^2 } G_{}(\vec{x}, \vec{y}) u_{n,0}(\vec{x})u_{n,i}(\vec{y}) d\vec{x} d\vec{y} \to 0,
\qquad i=0,1,2,
\end{align*}
by the Lebesgue dominated convergence theorem, since
$G_{}(\vec{x}, \vec{y}) \in L^1(\mathbb{T}^2\times \mathbb{T}^2;\mathbb{R}) $, $u_{n,i} $ is valued in $\{0,1\}$, $\|u_{n,0}\|_{L^1(\mathbb{T}^2)} \to 0$, thus $u_{n,0} \to 0$ a.e.. The interactions involving only
majority phases, i.e.
\begin{align*}
\int_{\mathbb{T}^2 \times \mathbb{T}^2 } G_{}(\vec{x}, \vec{y}) u_{n,i}(\vec{x})u_{n,j}(\vec{y}) d\vec{x} d\vec{y} ,
\qquad i,j=1,2,
\end{align*}
are estimated using
\ref{perimeter controls interaction - interaction estimate - with absolute values}:
summing over $i,j=1,2$ gives
\begin{align}
	\bigg|
\sum_{i,j=1}^{2} \gamma_{ij} \int_{\mathbb{T}^2 \times \mathbb{T}^2 } G_{}(\vec{x}, \vec{y}) u_{n,i}(\vec{x})u_{n,j}(\vec{y}) d\vec{x} d\vec{y} \bigg|
&\le
\sum_{i,j=1}^{2} |\gamma_{ij}| \bigg| \int_{\mathbb{T}^2 \times \mathbb{T}^2 } G_{}(\vec{x}, \vec{y}) u_{n,i}(\vec{x})u_{n,j}(\vec{y}) d\vec{x} d\vec{y}{\bigg|} \notag\\
& \le 
c_{3,R,\mathbb{T}^2}
\sum_{i,j=1}^{2}  | \gamma_{ij} | 
|\Omega_{n,i}| |\Omega_{n,j}|,
\label{compactness 0 order - interaction estimate}
\end{align}
where in the last inequality we applied \eqref{perimeter controls interaction - interaction estimate - with absolute values} with $E_1=\Omega_{n,i}$, $E_2=\Omega_{n,j}$, and the term $+\ln 3$ is due to  $\frac{1}{3}\le |\Omega_{n,1}|=|\Omega_{n,2}|$.
By \eqref{perimeter controls interaction - perimeter estimate}, the perimeter part is bounded from below by
\begin{align}
\label{compactness 0 order - perimeter estimate}
\frac{1}{2} (\Peri(\Omega_{n,1})+\Peri(\Omega_{n,2}))\ge \frac{c_{isoper, \mathbb{T}^2} }{2}  (\sqrt{ |\Omega_{n,1}|}  + \sqrt{ |\Omega_{n,2}|}),  
\end{align}
where the factor $1/2 $ is to account for double counting.
Since $\frac{1}{3}\le |\Omega_{n,1}|=|\Omega_{n,2}|\le\frac{1}{2}$,
\begin{align*}
\sum_{i,j=1}^{2} |\gamma_{ij}||\Omega_{n,i}| |\Omega_{n,j}|
&\le \frac{1}{2} [ (  | \gamma_{11} |  +  | \gamma_{12} | )
|\Omega_{n,1}|+ (  | \gamma_{12} | +   | \gamma_{22} |  )   |\Omega_{n,2}|]\\
&
\le \frac{1}{2\sqrt{2}} [ (|\gamma_{11}|+ |\gamma_{12}|)
\sqrt{|\Omega_{n,1}|}+ (|\gamma_{12}|+ |\gamma_{22}|) \sqrt{|\Omega_{n,2}|}]\\
&
\le \frac{|\gamma_{11}|+ |\gamma_{12}|}{2\sqrt{2} c_{isoper, \mathbb{T}^2}}  
\Peri(\Omega_{n,1})+ \frac{|\gamma_{12}|+ |\gamma_{22}|}{2\sqrt{2} c_{isoper, \mathbb{T}^2}}  \Peri(\Omega_{n,2}),
\end{align*}
hence \eqref{compactness 0 order - interaction estimate} becomes
\begin{align*}
\sum_{i,j=1}^{2} \gamma_{ij} &\int_{\mathbb{T}^2 \times \mathbb{T}^2 } G_{}(\vec{x}, \vec{y}) u_{n,i}(\vec{x})u_{n,j}(\vec{y}) d\vec{x} d\vec{y} 
\ge 
 - 
c_{3,R,\mathbb{T}^2}  \bigg[
\frac{|\gamma_{11}|+ |\gamma_{12}|}{2\sqrt{2} c_{isoper, \mathbb{T}^2}}  
\Peri(\Omega_{n,1})+ \frac{|\gamma_{12}|+ |\gamma_{22}|}{2\sqrt{2} c_{isoper, \mathbb{T}^2}}  \Peri(\Omega_{n,2})
\bigg],
\end{align*}
 and 
 combining with \eqref{compactness 0 order - perimeter estimate} gives
\begin{align*}
\frac{1}{2} \sum_{i=1}^{2} \Peri(  \Omega_{n,i}  )
&+ \sum_{i,j=1}^{2} \gamma_{ij} \int_{\mathbb{T}^2 \times \mathbb{T}^2 } G_{}(\vec{x}, \vec{y}) u_{n,i}(\vec{x})u_{n,j}(\vec{y}) d\vec{x} d\vec{y} 
\ge
\frac{1}{2}
\sum_{i=1}^{2} 
\bigg[ \underbrace{ 1-c_{3,R,\mathbb{T}^2}
  \frac{    |\gamma_{i1}|+ |\gamma_{i2}| }{ \sqrt{2} c_{isoper, \mathbb{T}^2}}  }_{>0 \text{ by \eqref{zeroth order gamma limit - negative interaction coefficients - conditions}}} \bigg]
\Peri(\Omega_{n,i}).
\end{align*}
Thus
uniform boundedness in energy implies uniform boundedness in $BV(\mathbb{T}^2)$
for $u_{\eta_{n,i}}$, $i=1,2$, which is thus (up to a subsequence which will not be relabeled) weakly converging in $BV(\mathbb{T}^2)$.

\medskip

$\Gamma-\liminf$ {\em inequality}. The perimeter part is clearly lower semicontinuous. We now show the interaction terms
are, up to subsequence, continuous with respect to weak convergence in $BV(\mathbb{T}^2)$. 
First, all interactions involving the minority phase are vanishing, as shown in the compactness part.
The interactions involving only majority phases are
\begin{align*}
\int_{\mathbb{T}^2 \times \mathbb{T}^2 } G_{}(\vec{x}, \vec{y}) u_{n,i}(\vec{x})u_{n,j}(\vec{y}) d\vec{x} d\vec{y} ,
\qquad i,j=1,2.
\end{align*}
The Green's function
$G_{}(\vec{x}, \vec{y}) \in L^1(\mathbb{T}^2\times \mathbb{T}^2;\mathbb{R})$, and $u_{n,i} \to u_{i}$ weakly in $BV(\mathbb{T}^2)$, and hence
strongly in $L^1(\mathbb{T}^2)$, and up to a subsequence (which we do not relabel) a.e.. Thus
\[G_{}(\vec{x}, \vec{y}) u_{n,i}(\vec{x})u_{n,j}(\vec{y}) \to G_{}(\vec{x}, \vec{y}) u_{i}(\vec{x})u_{j}(\vec{y}) \ a.e.,\qquad
G_{}(\vec{x}, \vec{y}) u_{n,i}(\vec{x})u_{n,j}(\vec{y}) \le G_{}(\vec{x}, \vec{y})  \in L^1(\mathbb{T}^2\times \mathbb{T}^2;\mathbb{R}),\]
hence by Lebesgue dominated convergence theorem we get 
\begin{align*}
\int_{\mathbb{T}^2 \times \mathbb{T}^2 } G_{}(\vec{x}, \vec{y}) u_{n,i}(\vec{x})u_{n,j}(\vec{y}) d\vec{x} d\vec{y}\to \int_{\mathbb{T}^2 \times \mathbb{T}^2 } G_{}(\vec{x}, \vec{y}) u_{i}(\vec{x})u_{j}(\vec{y}) d\vec{x} d\vec{y} ,
\qquad i,j=1,2,
\end{align*}
i.e. the interaction terms
are continuous with respect to weak convergence in $BV(\mathbb{T}^2)$. 
Now, we can infer the $\Gamma-\liminf$ inequality by contradiction: assume the opposite, i.e. there exists a sequence
$\eta_n$ such that
\[\liminf_{n\to +\infty} E_{\eta_n} (\Omega_{\eta_n,0},\Omega_{\eta_n,1},\Omega_{\eta_n,2})
< E_{L} (\Omega_{1},\Omega_{2}),\qquad  \Omega_{i}=\supp u_{i},\ \Omega_{\eta_n,i}=\supp u_{\eta_n,i},\quad i=0,1,2.\]
Then we can extract a subsequence $\eta_{n_k}$ such that
\[\liminf_{k\to +\infty} E_{\eta_{n_k}} (\Omega_{\eta_{n_k},0},\Omega_{\eta_{n_k},1},\Omega_{\eta_{n_k},2})
\ge E_{L} (\Omega_{1},\Omega_{2}),\]
since the perimeter part is clearly lower semicontinuous, and the interaction terms
are, up to subsequence, continuous with respect to weak convergence in $BV(\mathbb{T}^2)$. We have thus obtained 
a contradiction.

\medskip

$\Gamma-\limsup$ {\em inequality}. Let $u_i:\mathbb{T}^2\longrightarrow \{0,1\}$, $i=1,2$, be given,
satisfying 
$u_1+u_2=1 $.
 We want to construct 
a sequence of functions $u_{\eta,0},u_{\eta,1},u_{\eta,2}:\mathbb{T}^2\longrightarrow \{0,1\}$ such that
\[\|u_{\eta,0}\|_{L^1(\mathbb{T}^2)}=\eta^2 M,\quad \sum_{i=0}^{2}u_{\eta,i}=1.\]
Since $C^\infty(\mathbb{T}^2)$ is dense in $BV(\mathbb{T}^2)$, we can choose smooth functions $u_{\eta,1},u_{\eta,2} $
weakly converging to $u_1,u_2$, 
such that $\sum_{i=1}^{2}u_{\eta,i}=1$. Then, since $u_{\eta,1},u_{\eta,2} $ are smooth, we can easily
carve out a smooth $u_{\eta,0}:\mathbb{T}^2\longrightarrow \{0,1\}$ satisfying all the three mass constraints. Then
since by construction $|u_{\eta,i}|_{TV}\to |u_{i}|_{TV}$, $i=1,2$, and the interaction terms are continuous up to subsequence
with respect to the weak convergence in $BV(\mathbb{T}^2)$ as shown in the $\Gamma-\liminf$ part, we conclude that
the energies of $u_{\eta,0},u_{\eta,1},u_{\eta,2}$ converge to that of
$u_0,u_{1},u_{2}$, concluding the proof.
\end{proof}

 \begin{theorem}\label{1st order Gamma convergence - negative interaction coefficients}
 	The conclusion of Theorem \ref{1st order Gamma convergence} can be extended
 	to allow the interaction coefficients to be negative, as long as these satisfy
 	\begin{align}
 	\frac{ c_{3,R,\mathbb{T}^2} }{2} M\sum_{i=1}^2 |\Gamma_{i0}| + 2M^2|\Gamma_{00}|<s c_{isoper, \mathbb{T}^2} \sqrt{M} ,
 	\label{1st order Gamma convergence - negative interaction coefficients - conditions}
 	\end{align}
 	for some $s<1$, 
 	where $c_{3,R,\mathbb{T}^2}$ is defined in \eqref{zeroth order gamma limit - negative interaction coefficients - conditions}.
 \end{theorem}

\begin{proof}

	Consider a sequence $\eta_n\to 0$, and 
\[ u_{n,i}:\mathbb{T}^2 \longrightarrow \{0,1\},\ 	\sum_{i = 0}^2 u_{n,i} =1,
\ 	\Omega_{n,i}:=\supp u_{n,i},
\ |\Omega_{n,0}| = \eta_n^2 M,
\ \frac{1}{3}\le |\Omega_{n,1}|=|\Omega_{n,2}|=\frac{1-\eta_n^2M}{2},
\qquad i=0,1,2. \]
Set
\begin{align*}
E_\eta^{(1)}:=\frac{E_\eta -E_L}{\eta},
\end{align*}
i.e.
\begin{align*}
E_\eta^{(1)}(\Omega_0,\Omega_1,\Omega_2)=
\frac{\Peri(\Omega_0)}{\eta}
+ \sum_{i=1}^2  \frac{ \Gamma_{i0}}{\eta^2} \int_{\om_i} \int_{\om_0} G_{}(\vec{x}, \vec{y}) d\vec{x} d\vec{y}
+\frac{\Gamma_{00}}{\eta^4 |\ln \eta|} \int_{\om_0} \int_{\om_0} G_{}(\vec{x}, \vec{y}) d\vec{x} d\vec{y}.
\end{align*}

{\em Compactness.} We want to show that, if $\sup_n E_{\eta_n}^{(1)}(\Omega_{n,0},\Omega_{n,1},\Omega_{n,2})<+\infty$,
then there exist a Radon measures $\mu_i$ such that $\eta_n^{-2}u_{n,i}\to \mu_i$ in the weak-* topology. In order to rely
on previous results in literature, we need to ensure that the interaction terms can be controlled by the perimeters,
when the interaction coefficients are sufficiently small.

By the isoperimetric inequality, and $|\Omega_{n,0}|=\eta_n^2M $,
\begin{align}
\Peri(\Omega_{n,0}) & \ge c_{isoper, \mathbb{T}^2} \sqrt{|\Omega_{n,0}|} =c_{isoper, \mathbb{T}^2} \sqrt{M} \eta_n  .
\label{compactness 1 order - perimeter estimate}
\end{align}
%
%
Using \eqref{perimeter controls interaction - interaction estimate - with absolute values}
with $E_1=\Omega_{n,i}$, $E_2=\Omega_{n,0}$, and recalling that $|\Omega_{n,i}|\ge 1/3$,
we get
\begin{align} 
	\bigg|\frac{ \Gamma_{i0}}{\eta_n^2}\int_{\Omega_{n,0} \times \Omega_{n,i}} G_{}(\vec{x}, \vec{y}) d\vec{x} d\vec{y}\bigg|
\le \frac{ \Gamma_{i0}}{\eta_n^2}
c_{3,R,\mathbb{T}^2}
|\Omega_{n,i}||\Omega_{n,0}| \le \frac{M}{2}c_{3,R,\mathbb{T}^2}|\Gamma_{i0}|,\qquad
i=1,2.\label{compactness 1 order - interaction estimate 0-12}.
\end{align}
Similarly, using \eqref{perimeter controls interaction - interaction estimate - with absolute values}
with $E_1=E_2=\Omega_{n,0}$, and recalling that $|\Omega_{n,0}|= M\eta_n^2$,
we get
\begin{align}
\bigg|\int_{\Omega_{n,0} \times \Omega_{n,0}}  G_{}(\vec{x}, \vec{y}) d\vec{x} d\vec{y}\bigg|	
	&\le \bigg(\frac{1-\ln |\Omega_{n,0}|+\ln \pi}{4\pi }  +\sup_{\vec{x}, \vec{y}\in \mathbb{T}^2} |R(\vec{x}, \vec{y})|\bigg)
	|\Omega_{n,0}|^2\notag 
\end{align}
and recalling that $|\Omega_{n,0}|=\eta_n^2M \ll 1$, for all sufficiently small $\eta_n$,
such that $-\ln |\Omega_{n,0}| \ge 1+\ln \pi$,
 it holds
\begin{align}
	\bigg|	\int_{\Omega_{n,0} \times \Omega_{n,0}}  G_{}(\vec{x}, \vec{y}) d\vec{x} d\vec{y}\bigg|	&\le
	2|\ln \eta_n||\Omega_{n,0}|^2 \le 2|\ln \eta_n|\eta_n^4 M^2
	\label{compactness 1 order - interaction estimate 0-0}.
\end{align}
Combining with \eqref{compactness 1 order - perimeter estimate}, \eqref{compactness 1 order - interaction estimate 0-0},
\eqref{compactness 1 order - interaction estimate 0-12},
\begin{align*}
	\bigg|
\sum_{i=1}^2  \frac{ \Gamma_{i0}}{\eta_n^2}\int_{\Omega_{n,0} \times \Omega_{n,i}} G_{}(\vec{x}, \vec{y}) d\vec{x} d\vec{y}
&+\frac{\Gamma_{00}}{\eta_n^4 |\ln \eta_n|}\int_{\Omega_{n,0} \times \Omega_{n,0}} G_{}(\vec{x}, \vec{y}) d\vec{x} d\vec{y} \bigg|\\
&\le
\sum_{i=1}^2  \bigg| \frac{ \Gamma_{i0}}{\eta_n^2}\int_{\Omega_{n,0} \times \Omega_{n,i}} G_{}(\vec{x}, \vec{y}) d\vec{x} d\vec{y} \bigg| 
+\frac{|\Gamma_{00}|}{\eta_n^4 |\ln \eta_n|}	\bigg| \int_{\Omega_{n,0} \times \Omega_{n,0}} G_{}(\vec{x}, \vec{y}) d\vec{x} d\vec{y}	\bigg| 
\\
&\le
\frac{ c_{3,R,\mathbb{T}^2} }{2} M\sum_{i=1}^2 |\Gamma_{i0}| + 2M^2|\Gamma_{00}|
\overset{\eqref{1st order Gamma convergence - negative interaction coefficients - conditions}}{<} s c_{isoper, \mathbb{T}^2} \sqrt{M} 
\overset{\eqref{compactness 1 order - perimeter estimate}}{\le}  s \frac{\Peri(\Omega_{n,0})}{\eta_n} ,
\end{align*}
which in turn gives
\begin{align*}
E_{\eta_n}^{(1)}(\Omega_{n,0},\Omega_{n,1},\Omega_{n,2})
&=
\frac{\Peri(\Omega_{n,0})}{\eta_n} +
\sum_{i=1}^2  \frac{ \Gamma_{i0}}{\eta_n^2}\int_{\Omega_{n,0} \times \Omega_{n,i}} G_{}(\vec{x}, \vec{y}) d\vec{x} d\vec{y}
+\frac{\Gamma_{00}}{\eta_n^4 |\ln \eta_n|}\int_{\Omega_{n,0} \times \Omega_{n,0}} G_{}(\vec{x}, \vec{y}) d\vec{x} d\vec{y}
\\
&\ge
\frac{\Peri(\Omega_{n,0})}{\eta_n} -
\bigg| \sum_{i=1}^2  \frac{ |\Gamma_{i0}|}{\eta_n^2}\int_{\Omega_{n,0} \times \Omega_{n,i}} G_{}(\vec{x}, \vec{y}) d\vec{x} d\vec{y}
+\frac{|\Gamma_{00}|}{\eta_n^4 |\ln \eta_n|}\int_{\Omega_{n,0} \times \Omega_{n,0}} G_{}(\vec{x}, \vec{y}) d\vec{x} d\vec{y} \bigg| 
\\
&\ge
(1-s)\frac{\Peri(\Omega_{n,0})}{\eta_n},
\end{align*}
which enable us to use the results from \cite[Sections~4--6]{bi1}, to infer compactness.
The proof of $\Gamma-\liminf$ (resp. $\Gamma-\limsup$) inequality then follows the same arguments from
\cite[Theorem~6.1]{bi1} (resp. Theorem \ref{1st order Gamma convergence}).

\end{proof}

\end{document}